\documentclass[11pt,leqno]{article}

\usepackage{amsmath,amsthm,amsfonts}
\usepackage {latexsym}
\usepackage{amssymb}

\newcommand{\diam}{\mbox{\rm diam}}

\newcommand{\bse}{\begin{equation}}

\renewcommand{\mod}{{\rm mod\ }}

\newcommand{\bea}{\begin{eqnarray}}
\newcommand{\eea}{\end{eqnarray}}
\newcommand{\be}{\begin{equation}}
\newcommand{\ee}{\end{equation}}

\newcommand{\half}{\frac{1}{2}}

\newcommand{\eps}{{\varepsilon}}
\newcommand{\R}{{\mathbb R}}
\newcommand{\tor}{{\mathbb T}}
\newcommand{\Z}{{\mathbb Z}}

\newcommand{\Compl}{{\mathbb C}}
\newcommand{\Erw}{\,{\mathbb E}\,}
\newcommand{\Prob}{{\mathbb P}}

\newcommand{\calP}{{\mathcal P}}
\newcommand{\calF}{{\mathcal F}}
\newcommand{\calS}{{\mathcal S}}

\newcommand{\calE}{{\mathcal E}}

\newcommand{\gtr}{\gtrsim}

\newtheorem{theorem}{Theorem}[section]
\newtheorem{lemma}[theorem]{Lemma}
\newtheorem{cor}[theorem]{Corollary}
\newtheorem{prop}[theorem]{Proposition}

\theoremstyle{definition}
\newtheorem{defi}[theorem]{Definition}

\theoremstyle{remark}
\newtheorem{remark}[theorem]{Remark}

\def\summe{\sum_{n=1}^N}
\def\les{\lesssim}
\def\gtr{\gtrsim}
\def\nn{\nonumber}
\def\la{\langle}
\def\ra{\rangle}
\def\Impl{\;\Longrightarrow\;}
\def\calD{{\mathcal D}}
\def\what{\widehat}

\numberwithin{equation}{section}

\begin{document}

\title{On the Hardy-Littlewood majorant problem for random sets}
\author{G.\ Mockenhaupt, W.\ Schlag}
\date{}
\maketitle

\section{The  majorant property: Some generalities}

This paper is concerned with
the majorant property of various randomly generated subsets of $[1,N]$.
More precisely, suppose $A_N\subset [1,N]$ is a sequence of sets so that
$|A_N|\asymp N^{\rho}$ for some fixed~$0<\rho<1$ as $N\to\infty$. 
For example, one can take $A_N$ to be the squares, cubes, etc., or
(multi-dimensional) arithmetic progressions. 
As in~\cite{M}, given $p\ge 2$, one asks for the smallest power~$\gamma=\gamma(p)>0$ 
(which might be also specific to the sequence~$A_N$) 
such that
\begin{equation}
\label{eq:basic_ques}
\sup_{|a_n|\le 1} \Big\|\sum_{n\in A_N} a_n e(n\cdot) \Bigr\|_p
\le C\, N^\gamma \Big\|\sum_{n\in A_N}  e(n\cdot) \Bigr\|_p. 
\end{equation}
This is only one out of several ways of stating the majorant problem.
\cite{M}~asks for a power $\gamma$ that applies to \underline{all} 
$A_N\subset [1,N]$ simultaneously. 
If $p$ is an even integer, then one can take $C=1$ and $\gamma=0$
as realized by Hardy and Littlewood. 
On the other hand, it has also been known for some time that one cannot take~$\gamma=0$
if~$p$ is not an even integer. Moreover, a quantitative lower bound 
of~$\exp(c\frac{\log N}{\log\log N})$ is obtained in~\cite{M} for~\eqref{eq:basic_ques}
with a particular choice of~$A_N$. 
If~\eqref{eq:basic_ques} holds for all~$\gamma>0$ and appropriate~$p$, then it would
imply the restriction and therefore also the Kakeya conjecture, see~\cite{M}
 for those matters.
One always has the bound
\[ 
\sup_{|a_n|\le 1} \Big\|\sum_{n\in A_N} a_n e(n\cdot) \Bigr\|_p 
\le C \Bigl(\frac{N}{|A_N|}\Bigr)^{\frac{1}{p}} \Big\|\sum_{n\in A_N}  e(n\cdot) \Bigr\|_p
\]
by Hausdorff-Young and the obvious lower bound 
$\Big\|\sum_{n\in A_N}  e(n\cdot) \Bigr\|_p^p\gtrsim |A_N|^p N^{-1}$. 
This settles the case of any sequence of large sets, i.e., $\rho=1$, as well as 
all arithmetic progressions. Another easy estimate can be obtained by interpolation.
Indeed, if $2<p<4$, say, then interpolating between $2$ and~$4$ yields 
$\gamma \le (1-\frac{p}{4})(1-\frac{2}{p})$. 
It turns out that this interpolation can be done more carefully, which gives optimal
results for sets $A_N$ whose Dirichlet kernel satisfies a certain ``reverse interpolation inequality.''
To this end, let 
$ \calP_A:= \{ \sum_{n\in A} a_n e(n\theta)\:|\: |a_n| \le 1 \}. $
Then, with $A=A_N$ for simplicity, for any odd integer $p>2$, 
\bea
&& \sup_{|a_n|\le 1} \int_0^1 \Bigl|\sum_{n\in A} a_n e(n\theta) \Bigr|^p\,d \theta
= \sup_{|a_n|\le1} \sum_{n\in A} a_n \int_0^1 e(n\theta) \sum_{k\in A} \bar{a}_k e(-k\theta) 
\Big| \sum_{\ell\in A} a_\ell e(\ell\theta)\Big|^{p-2}\,d\theta \nn \\
&& \le \sup_{g\in \calP_A} \; \sqrt{|A|}\,\left( \sum_{n\in A} \big| \what{g|g|^{p-2}} (n) \big|^2 \right)^{\half} \le  \sup_{g\in \calP_A} \; \sqrt{|A|}\,\| g\|_{2(p-1)}^{p-1} \label{eq:planch} \\
&& \le \Big\|\sum_{n\in A} e(n\cdot) \Big\|_2 \Big\|\sum_{n\in A} e(n\cdot) \Big\|_{2(p-1)}^{p-1} 
\label{eq:maj_step}.
\end{eqnarray}
Here the first inequality sign in~\eqref{eq:planch} follows by putting absolute values
inside and Cauchy-Schwarz, the second is Plancherel, and~\eqref{eq:maj_step} uses
the majorant property on~$2(p-1)$. Now assume the following condition
\begin{equation}
\label{eq:star}
\Big\|\sum_{n\in A_N} e(n\cdot) \Big\|_2 \Big\|\sum_{n\in A_N} e(n\cdot) \Big\|_{2(p-1)}^{p-1} 
\le C_\eps \, N^{\eps} \Big\|\sum_{n\in A_N} e(n\cdot) \Big\|_p^p.
\end{equation}
In view of the preceding, one then has~\eqref{eq:basic_ques} for any $\gamma>0$. 
This condition, which is of basic importance for most of our work, is basically the
reverse of the usual interpolation inequality. 
One checks immediately that arithmetic progressions satisfy~\eqref{eq:star}.
Also, observe that
any sequence $A_N$ for which~\eqref{eq:star} holds for all~$p$ satisfies~\eqref{eq:basic_ques}
for all $p$ with~$\gamma>0$. Indeed, this follows inductively from the argument leading
up to~\eqref{eq:maj_step} using the majorant property from the previous stage $2(p-1)$  to pass
to the next stage~$p$. Finally, interpolation is required to obtain the desired bound
for all~$p$ (at the cost of~$N^\eps$). 
Another case which is covered by this argument, but not the previous one based on Hausdorff-Young, are
multi-dimensional arithmetic progressions. For example, one easily checks that 
\begin{equation}
\label{eq:multi} 
A=\{b+j_1 a_1+ j_2 a_2\:|\: 0\le j_1< L_1,\,0\le j_2< L_2\}
\end{equation}
with $a_1L_1 <a_2$, satisfies 
\[ \Big\|\sum_{n\in A} e(n\cdot) \Big\|_p^p \asymp (L_1 L_2)^{p-1} \]
for $p>1$. 
Another interesting case are the squares $A_N=\{n^2\:|\: 1\le n\le \sqrt{N}\}$. 
In this case it is well-known that the there is a ``kink'' at $p=4$, 
\bea
 \Big\|\sum_{n\in A_N} e(n\cdot) \Big\|_p &\le& C_\eps\,N^{\eps+\half} \text{\ \ if\ }2\le p\le 4, \nn \\
 \Big\|\sum_{n\in A_N} e(n\cdot) \Big\|_p &\le& C_\eps\,N^{1-\frac{2}{p}+\eps} \text{\ \ if\ }p\ge 4, \nn
\end{eqnarray}
so that \eqref{eq:star} holds only for $2\le p\le 3$. In particular, the argument leading
up to~\eqref{eq:maj_step} gives the (trivial) statement that the majorant property
holds at~$p=3$ for the squares. A nontrivial statement can be obtained by improving on
the use of Plancherel in~\eqref{eq:planch}. Indeed, it is a well-known fact that
\begin{equation}
\label{eq:dual_4} 
\Big\| \sum_{n=1}^N a_n e(n^2\theta) \Big\|_4 \le C_\eps\,N^\eps\, \left(\sum_{n=1}^N |a_n|^2 \right)^{\half} \Longleftrightarrow 
\left(\sum_{n=1}^N |\hat{f}(n^2)|^2 \right)^{\half} \le C_\eps\,N^\eps \|f\|_{L^{\frac43}(\tor)},
\end{equation}
the second statement being the dual of the first. 
This can be checked by reducing the $L^4$-norm to an $L^2$-norm by squaring, and then using
Cauchy-Schwarz and the $N^\eps$-bound on the divisor function, see~\cite{B3}. 
We now repeat the argument leading up to~\eqref{eq:maj_step} to conclude the following.
Let
\[\calP:= \{ \summe a_n e(n^2\theta)\:|\: |a_n| \le 1 \}. \]
If $p=3k+1$, then one can apply the majorant property at $\frac{4}{3}(p-1)$ so that 
\bea
&& \sup_{|a_n|\le 1} \int_0^1 \Bigl|\summe a_n e(n^2\theta) \Bigr|^p\,d \theta
= \sup_{|a_n|\le1} \summe a_n \int_0^1 e(n^2\theta) \sum_{k=1}^N \bar{a}_k e(-k^2\theta) 
\Big| \sum_{\ell=1}^N a_\ell e(\ell^2\theta)\Big|^{p-2} \, d\theta\nn \\
&& \le \sup_{g\in \calP} \; \sqrt{|A|}\;\left( \summe \big| \what{g|g|^{p-2}} (n^2) \big|^2 \right)^{\half} \le  \sup_{g\in \calP} \; \sqrt{|A|}\;\| g\|_{\frac{4}{3}(p-1)}^{p-1} \label{eq:quadrat} \\
&& \le 
\Big\|\summe e(n^2\cdot) \Big\|_2 \Big\|\summe e(n^2\cdot) \Big\|_{\frac{4}{3}(p-1)}^{p-1}
\le C_\eps\, N^\eps N^{\half} N^{p-\frac52}\le C_\eps\,N^\eps N^{p-2}\nn \\
&& \le C_\eps\,N^\eps\,\Big\|\summe e(n^2\cdot) \Big\|_p^p . \nn
\end{eqnarray}
Here we used~\eqref{eq:dual_4} in~\eqref{eq:quadrat}. 
This implies that for the sequence of squares~\eqref{eq:basic_ques} holds with any $\gamma>0$ 
at $p=7,13,19$ etc.

Another case of sets $A_N$ that \underline{do not} satisfy~\eqref{eq:star} are
random subsets~$A_N\subset [1,N]$. Indeed, we show below that random sets $A_N$ which are
obtained by selecting each integer $1\le n\le N$ with probability~$\tau$ have the property
that for $p>1$ 
\[ \Erw \Big\|\sum_{n\in A_N} e(n\theta) \Big\|_p^p \asymp \tau^p N^{p-1} + (\tau N)^{\frac{p}{2}}, \]
see Theorem~\ref{thm:Lp}. 
 The two terms on the right balance at $\tau_{crit}=N^{-1+\frac{2}{p}}$ so that 
it is clear that~\eqref{eq:star} cannot hold in general. 
The main objective of the following section is to show that nevertheless, 
such random subsets \underline{do satisfy}~\eqref{eq:basic_ques} with large probability. 
The method to some extent resembles the calculation from~\eqref{eq:maj_step}, but
is of course more involved. We rely on a probabilistic lemma from Bourgain's work~\cite{B1}. 

It is possible to abstract the arguments below, and then verify that
various examples satisfy the conditions of such an abstract theorem, 
the most important one being condition~\eqref{eq:star}. More precisely, 
starting with a {\em deterministic} sequence $A_N$, 
define $\calS_N(\omega)=\{n\in A_N\:|\: \xi_n=1\}$ where $\xi_n$ are i.i.d.~selector variables 
satisfying 
$\Prob[\xi_n=1]=\tau=1-\Prob[\xi_n=0]$. If, amongst other things, \eqref{eq:star} 
holds for~$A_N$, then much of what is done 
in the following section goes through. On the other hand, some improvements which
we obtain below for the case of arithmetic progressions are not easily axiomatized. 
Moreover, since we do not have any examples apart from (multi-dimensional) arithmetic
progressions, we have decided against casting this into a more general framework.
Thus, we write out the main argument only for arithmetic progressions.
If~\eqref{eq:star} is violated, then our method applies only to certain~$p$ or
after suitable modifications. For example, one can check that 
the machinery which we develop below shows that with high probability 
random subset of the squares satisfy~\eqref{eq:basic_ques} at~$p=7$ for any $\gamma>0$.
 This requires invoking the (almost) $\Lambda(4)$ property of the squares as 
in~\eqref{eq:quadrat}. It seems difficult to obtain the desired bound for all~$p$ 
in case of the squares. 

In addition to random subsets we also consider perturbations of arithmetic
progressions. This means that each element of a given arithmetic progression is
shifted independently and randomly by some small amount. We again show that most sets
obtained in this fashion satisfy~\eqref{eq:basic_ques} for any $\gamma>0$, see Theorem~\ref{thm:AP1}. 
As before, the method can be presented abstractly for perturbations
of arbitrary sets~$A_N$ that satisfy condition~\eqref{eq:star}.

\section{Random subsets have the majorant property} 
\label{sec:RMP1}

\begin{theorem}
\label{thm:Lp} 
Let $0<\delta<1$ be fixed. For every positive integer $N$ we let $\xi_j=\xi_j(\omega)$ be i.i.d.~variables with $\Prob[\xi_j=1]=\tau$, $\Prob[\xi_j=0]=1-\tau$ where $\tau=N^{-\delta}$. Define a random subset
\[ S(\omega)=\{j\in[1,N]\:|\: \xi_j(\omega)=1\}.\]
Then for every $\eps>0$ and $7\ge p\ge2$  one has 
\begin{equation}
\label{eq:prob} 
\Prob\Bigl[\sup_{|a_n|\le1} \Big\|\sum_{n\in S(\omega)} a_n e(n\theta) \Big\|_{L^p(\tor)} \ge N^\eps\,
             \Big\|\sum_{n\in S(\omega)}  e(n\theta) \Big\|_{L^p(\tor)} \Bigr] \to 0
\end{equation}
as $N\to\infty$. Moreover, under the additional 
restriction~$\delta\le\half$, \eqref{eq:prob} holds for all $p\ge 7$.
\end{theorem}

\noindent We show below that the $N^\eps$-factor can be removed in certain cases, for example
when~$p=3$. The restriction $\delta\le \half$ for $p\ge 7$ 
appears to be of a purely technical nature, and we believe that the theorem should hold for 
all~$\delta\in(0,1)$. 

The proof of Theorem~\ref{thm:Lp} relies on a method that Bourgain developed for the $\Lambda(p)$ 
problem, see~\cite{B1} and~\cite{B2}. 
In fact, in this situation we can avoid several complications that arose
in Bourgain's work. 
Notice that our Theorem~\ref{thm:Lp} is implied by Bourgain's existence theorem
of $\Lambda(p)$ sets provided $\delta\ge 1-\frac{2}{p}$, but not for $\delta<1-\frac{2}{p}$. 
Indeed, in the former case the random set $S$ will typically have cardinality $N^{\frac{2}{p}}$ 
or smaller, and such sets were shown by Bourgain~\cite{B1} to be $\Lambda(p)$-sets with large probability. 

\subsection{Random sums over asymmetric Bernoulli variables}

We first dispense with some simple technical statements about the behavior of
random sums with asymmetric Bernoulli variables as summands.
They are definitely standard, but lacking a precise reference we prefer to present them.

\begin{lemma}
\label{lem:kinchin}
Let $\eta_j$ be i.i.d.~variables so that $\Prob[\eta_j=1-\tau]=\tau$, $\Prob[\eta_j=-\tau]=1-\tau$.
Here $0<\tau<1$ is arbitrary. Let  $N\ge 1$ and
$\{a_j\}_{j=1}^N \in\Compl$ be given. Define $\sigma^2 = \tau(1-\tau) \sum_{j=1}^N |a_j|^2$. 
Then for $\lambda>0$,
\[ \Prob\Big[ \Big|\sum_{j=1}^N a_j\eta_j\Big| > \lambda\sigma \Big] \le 4e^{-\frac{\lambda^2}{8}}\]
provided 
\begin{equation}
\label{eq:LDTcond}
\max_{1\le j\le N} {\lambda|a_j|} \le 4\sigma.
\end{equation}
\end{lemma}
\begin{proof}
Assume first that all $a_j\in\R$. Then for any $t>0$
\bea 
&& \Prob\Big[\sum_{j=1}^N a_j\eta_j  > \lambda\sigma \Big] \le e^{-t\lambda\sigma} \Erw 
\exp\Big(t\sum_{j=1}^N a_j\eta_j \Big) \label{eq:pos_first} \\
&& = e^{-t\lambda\sigma} \prod_{j=1}^N \Big[\tau e^{(1-\tau)t a_j} + (1-\tau) e^{-\tau t a_j} \Bigr] 
\label{eq:exp_erw}
\end{eqnarray}
Next, we claim that 
\begin{equation}
\label{eq:square}
\tau e^{(1-\tau)x} + (1-\tau) e^{-\tau x} \le \exp(2\tau(1-\tau)x^2) \text{\ \ for all\ \ }|x|\le1.
\end{equation}
Observe that this property fails for $x=\tau^{-\half}$. To prove this, set
\[ 
\phi_\tau(x)=\exp(2\tau(1-\tau)x^2) - \tau e^{(1-\tau)x} - (1-\tau) e^{-\tau x}.
\]
By symmetry it suffices to consider the case $0\le x\le 1$ and to show that $\phi_\tau\ge0$ there.
Clearly,
\bea 
\phi_\tau'(x) &=& \tau(1-\tau)\big[4x\exp(2\tau(1-\tau)x^2) -  e^{(1-\tau)x} + e^{-\tau x}\big ] \nn \\
&\ge& \tau(1-\tau)\big[4x - e^{(1-\tau)x} +  e^{-\tau x}\big ] \label{eq:diff}
\end{eqnarray}
Differentiating the expression in brackets yields 
\[
4-(1-\tau)e^{(1-\tau)x}-\tau e^{-\tau x}\ge 4-(1-\tau)e^{(1-\tau)x}-\tau e^{(1-\tau) x}\ge 4 - e > 0
\]
for all $0\le x\le1$.  It follows that $\phi_\tau'(x)\ge\phi_\tau'(0)=0$ for $0\le x\le1$. 
Hence also $\phi_\tau(x)\ge\phi_\tau(0)=0$ for~$0\le x\le1$, as desired.
Inserting~\eqref{eq:square} into~\eqref{eq:exp_erw} gives
\[
 \Prob\Big[\sum_{j=1}^N a_j\eta_j  > \lambda\sigma \Big] \le \min_{t>0} e^{-t\lambda\sigma} 
\exp(2t^2 \sigma^2) = e^{-\frac{\lambda^2}{8}}
\]
provided for the minimizing choice of $t=t_0$ one has $\max_{j} |t_0\,a_j|\le 1$. But 
$t_0=\frac{\lambda}{4\sigma}$ and this condition therefore reads
\[ \max_{1\le j\le N} \frac{|\lambda||a_j|}{4\sigma} \le 1,\]
which is precisely \eqref{eq:LDTcond}.
Evidently, the same bound also holds for deviations less than $-\lambda\sigma$, which gives
$2e^{-\lambda^2/8}$ as an upper bound on the large deviation probability in the real case. Finally,
if $a_n\in\Compl$, then one splits into real and complex parts.
\end{proof}
 
\noindent 
Lemma~\ref{lem:kinchin} immediately leads to the following version of the Salem--Zygmund inequality for
asymmetric variables.

\begin{cor}
\label{cor:salem}
With $\eta_n$ and $\sigma$ as in the previous lemma 
\[ 
\Prob\Big[ \sup_{\theta\in\tor} \Big|\sum_{n=1}^N a_n\,\eta_n\, e(n\theta) \Big| > 20\,\sigma\sqrt{\log N}\Big] \le 4N^{-8}
\]
for any $a_n\in\Compl$ provided the following conditions hold:
\bea
 \sup_{1\le n\le N} 10|a_n|^2\log N &\le& \sigma^2 = \tau(1-\tau)\sum_{k=1}^N |a_k|^2 \nn \\
 10 &\le& \tau(1-\tau)N\log N \label{eq:salem_cond}. 
\end{eqnarray}
\end{cor}
\begin{proof} 
Let $\{\theta_j\}_{j=1}^{N^2}\subset\tor$ be a $N^{-2}$-net. Denote
\[ T_{N,\omega}(\theta) := \sum_{n=1}^N a_n\,\eta_n(\omega)\, e(n\theta).\]
Then, by Bernstein's inequality, and with the usual de la Vallee-Poussin kernel $V_N$, 
\bea
&& \min_{j} |T_{N,\omega}(\theta)-T_{N,\omega}(\theta_j)| \le N^{-2} \| T_{N,\omega}'\|_\infty \nn \\
&& \le N^{-2} \| T_{N,\omega}\|_2 \|V_N'\|_2 \le N^{-2} \Big(\sum_{n=1}^N|a_n|^2\Big)^{\half}\; 8\pi N^{\frac32} \nn \\
&& = \frac{8\pi \sigma N^{-\frac12}}{\sqrt{\tau(1-\tau)}} \le 10\sigma\sqrt{\log N}.\nn
\end{eqnarray}
The final inequality here follows from our assumption~\eqref{eq:salem_cond}.
Therefore, by Lemma~\ref{lem:kinchin},
\bea
\Prob\Big[ \sup_{\theta\in\tor} \Big|\sum_{n=1}^N a_n\,\eta_n\, e(n\theta) \Big| > 20\,\sigma\sqrt{\log N}\Big]  &\le& \sum_{j=1}^{N^2} \Prob\Big[\Big|\sum_{n=1}^N a_n\,\eta_n\, e(n\theta_j) \Big| > 10\,\sigma\sqrt{\log N}\Big] \nn \\
&\le& 4 N^2 \exp(-100\log N/8)\le 4 N^{-8}, \nn
\end{eqnarray}
which is precisely the bound claimed in the lemma. 
The first condition in~\eqref{eq:salem_cond} ensures that~\eqref{eq:LDTcond} holds.
\end{proof}

\noindent
In the proof of Theorem~\ref{thm:Lp} we shall need to know 
the typical size of the easier norm in~\eqref{eq:prob}. We determine this
norm in the following lemma.

\begin{lemma}
\label{lem:norm1} 
Let $\xi_j$ be selector variables as above with $\tau=N^{-\delta}$, $0<\delta<1$ fixed.
Let $p\ge2$ and define
\[ I_{p,N}(\omega) = \int_0^1 \Bigl|\sum_{n=1}^N \xi_n(\omega) e(n\theta)\Bigr|^p \, d\theta.\]
Then for some constants $C_p$, 
\[ 
C_p^{-1}\,\Bigl( \tau^p N^{p-1} + (\tau N)^{\frac{p}{2}}\Bigr) \le \Erw I_{p,N} \le C_p\, \Bigl(\tau^p N^{p-1} + (\tau N)^{\frac{p}{2}}\Bigr).
\]
Moreover,  there is some small constant $c_p$ such that
\[ \Prob\Big[I_{p,N} \le c_p( \tau^p N^{p-1} + (\tau N)^{\frac{p}{2}})\Big] \to 0\]
as $N\to\infty$. 
\end{lemma}
\begin{proof} Let $\eta_n(\omega)=\xi_n(\omega)-\tau$, so that $\Erw \eta_n=0$ and
$\Erw\eta_n^2 = \tau(1-\tau)$. Then
\bea
I_{p,N}(\omega) &\les& \int_0^1 \Big| \summe \tau e(n\theta)\Big|^p \,d\theta + 
\int_0^1 \Big| \summe \eta_n e(n\theta)\Big|^p \,d\theta  \nn \\
&\les& \tau^p N^{p-1} + \int_0^1 \Big| \summe \eta_n e(n\theta)\Big|^p \,d\theta. \label{eq:obere}
\end{eqnarray}
One now checks that 
\[ \Erw \int_0^1 \Bigl|\sum_{n=1}^N \eta_n e(n\theta)\Bigr|^p \, d\theta \le C_p\,(N\tau(1-\tau))^{\frac{p}{2}}.
\]
This can be verified by expanding the norm for even $p$ and then interpolating.
Indeed,
\bea
&& \Erw \int_0^1 \Big|\sum_{n=1}^N \eta_n e(n\theta)\Bigr|^{2k} \, d\theta =
\Erw \int_0^1 \Big|\sum_{n_1,\ldots,n_k=1}^N \eta_{n_1}\ldots\eta_{n_k} 
e((n_1+\ldots+n_k)\theta)\Bigr|^{2} \, d\theta \nn \\
&& = \sum_n \Erw \Bigl| \sum_{n_1+\ldots+n_k=n} \eta_{n_1}\ldots\eta_{n_k} \Bigr|^2
= \sum_{n_1 +\ldots +n_k=m_1+\ldots+m_k} \Erw [\eta_{n_1}\ldots\eta_{n_k}\eta_{m_1}\ldots\eta_{m_k}] \nn \\
&& \le C_k \sum_{r=1}^k \sum_{\substack{n_1,\ldots,n_r=1\\ s_1+\ldots+s_r=2k,\;s_i\ge 2}}^N
\Erw |\eta_{n_1}|^{s_1}\cdot\ldots\cdot \Erw|\eta_{n_r}|^{s_r} \label{eq:comb_fac1} \\
&&\le C_k \sum_{r=1}^k N^r(\tau(1-\tau))^r \le C_k (N\tau(1-\tau))^k. \label{eq:comb_fac2}
\end{eqnarray} 
The constants in \eqref{eq:comb_fac1} and~\eqref{eq:comb_fac2} are of a combinatorial nature and
not necessarily the same. The relevant point in~\eqref{eq:comb_fac1} is that $s_i\ge2$ which
is due to independence and $\Erw \eta_j=0$. In particular, $s_i\ge2$ implies the important fact $r\le k$.
Moreover, to pass to the last line we used that for every positive integer $s\ge2$ 
\[ \tau(1-\tau)\ge \Erw \eta_j^{s} = \tau(1-\tau)(\tau^{s-1}+(1-\tau)^{s-1}) 
    \ge 2^{2-s}\tau(1-\tau).
\]
To obtain the lower bound on the expectation, one splits the integral in $\theta$ into the region
where the Dirichlet kernel dominates the mean zero random sum and vice versa. More precisely,
with $h=\sqrt{\tau N^{-1}}=N^{-\frac{1+\delta}{2}}$, 
\bea
I_{p,N} &\gtr& \int_{|\theta|<\frac1N} \Big|\summe \tau e(n\theta) \Big|^p \,d\theta 
- \int_{|\theta|<\frac1N} \Big|\summe \eta_n e(n\theta) \Big|^p \,d\theta \nn \\
&& + \int_{h}^{1-h} \Big|\summe \eta_n e(n\theta) \Big|^p \,d\theta - \int_{h}^{1-h} \Big|\summe \tau e(n\theta) \Big|^p \,d\theta \nn \\
&\gtr& \tau^p N^{p-1} 
- C \int_{|\theta|<\frac1N} \Big|\summe \eta_n e(n\theta) \Big|^p \,d\theta  + \int\limits_{|\theta|>h} \Big|\summe \eta_n e(n\theta) \Big|^p \,d\theta - C\tau^p\,h^{1-p}. \label{eq:untere}
\end{eqnarray}
According to Corollary~\ref{cor:salem}, the first integral in \eqref{eq:untere} is 
\begin{equation}
\label{eq:fehler1} 
\les N^{-1}(\log N)^{\frac{p}{2}} (\tau(1-\tau) N)^{\frac{p}{2}}
\end{equation}
up to a negligible probability. For the second, one has because of $p\ge2$
\bea
\int_{h}^{1-h} \Big|\summe \eta_n e(n\theta) \Big|^p \,d\theta &\ge& \int_{0}^{1} \Big|\summe \eta_n e(n\theta) \Big|^p \,d\theta - \int\limits_{|\theta|\le h} \Big|\summe \eta_n e(n\theta) \Big|^p \,d\theta \nn \\
&\ge& \left(\int_{0}^{1} \Big|\summe \eta_n e(n\theta) \Big|^2 \,d\theta\right)^{\frac{p}{2}} - 
\int\limits_{|\theta|\le h} \Big|\summe \eta_n e(n\theta) \Big|^p \,d\theta \nn \\
&\ge& \left( \sum_{n=1}^N \eta_n^2 \right)^{\frac{p}{2}} - C\,h\,(N\tau\,\log N)^{\frac{p}{2}} 
\label{eq:split_2int}
\end{eqnarray}
where the last term in~\eqref{eq:split_2int} is obtained from Corollary~\ref{cor:salem}.
Using  $p\ge2$ again, 
\[
 \Erw \left( \sum_{n=1}^N \eta_n^2 \right)^{\frac{p}{2}} \ge \left( \Erw \sum_{n=1}^N \eta_n^2 \right)^{\frac{p}{2}} 
 \ge \Big(N\tau(1-\tau)\Big)^{\frac{p}{2}}. 
\]
In fact, Lemma~\ref{lem:kinchin} gives the following more precise estimate:
\begin{equation}
\Prob\Big[ \Big|\sum_{n=1}^N (\eta_n^2-\Erw \eta^2_n)\Big| \ge
\lambda \sqrt{N\Erw(|\eta_1^2-\Erw\eta_1^2|^2)}  \Big]\le 4e^{-\lambda^2/8}
\label{eq:LDT2}
\end{equation}
provided the conditions \eqref{eq:LDTcond} hold. One checks that
$\Erw( |\eta_1^2-\Erw \eta_1^2|^2) \asymp \tau(1-\tau)$. Hence it follows from~\eqref{eq:LDT2} that
 for large~$N$
\bea
&& \Prob \Big[ \sum_{n=1}^N \eta_n^2 \le \half \Erw \sum_{n=1}^N \eta_n^2 = \half N\tau(1-\tau)\Big] 
 \le \Prob \Big[ \Big|\sum_{n=1}^N \eta_n^2 -\Erw \sum_{n=1}^N \eta_n^2\Big|\ge \half N\tau(1-\tau)\Big] \nn \\
&& \le \Prob \Big[ \Big|\sum_{n=1}^N \eta_n^2 -\Erw \sum_{n=1}^N \eta_n^2\Big|\ge \log N \sqrt{N\tau(1-\tau)}\Big] 
 \le 4e^{-(\log N)^2/8}, \nn
\end{eqnarray}
since with our choice of parameters \eqref{eq:LDTcond} hold for large~$N$.
Inserting this bound into \eqref{eq:split_2int} now yields (recall that $h=\sqrt{\tau N^{-1}}=N^{-\frac{1+\delta}{2}}$)
\[ 
\int_{h}^{1-h} \Big|\summe \eta_n e(n\theta) \Big|^p \,d\theta \ge \left(\half N\tau(1-\tau)\right)^{\frac{p}{2}} - C\,N^{-\frac{1+\delta}{2}}\, (N\tau\,\log N)^{\frac{p}{2}} \gtrsim (N\tau)^{\frac{p}{2}}
\]
up to negligible probability. 
In view of this bound and \eqref{eq:fehler1}, one obtains from~\eqref{eq:untere} that
\bea
I_{p,N} &\gtrsim & 
\tau^p N^{p-1} 
- C \int_{|\theta|<\frac1N} \Big|\summe \eta_n e(n\theta) \Big|^p \,d\theta  + \int\limits_{|\theta|>h} \Big|\summe \eta_n e(n\theta) \Big|^p \,d\theta - C\tau^p\,h^{1-p} \nn \\
&\gtr& \tau^pN^{p-1} + (N\tau)^{\frac{p}{2}} \nn
\end{eqnarray}
up to negligible probability. To remove the final term in the first line we used that $(N\tau)^{\frac{p}{2}} \gg \tau^p h^{1-p}$
which follows from our choice of $h$ provided $N$ is big.
\end{proof}

\subsection{Suprema of  random processes}

We now collect the statements from Bourgain's paper that we will need.
The first is Lemma~1 from~\cite{B1} with $q_0=1$. 
In fact, Bourgain's lemma is slightly stronger because of certain $\log \frac{1}{\tau}$-factors.
While these factors are important for his purposes, they play no role in our argument.
We present the proof for the reader's convenience, following Bourgain's original argument.
Another proof was found by Ledoux and Talagrand~\cite{LT} which is close to the ideology
surrounding Dudley's theorem on suprema of Gaussian processes. While their point of view is perhaps more
conceptual, we have found it advantageous to follow~\cite{B1}.
Throughout, if $x\in\R^N$, then $|x|=|x|_{\ell^2_N}=\Big(\sum_{j=1}^N x_j^2\Big)^{\half}$
is the Euclidean norm. Secondly, $N_2(\calE,t)$ refers to the $L^2$-entropy
of the set~$\calE$ at scale~$t$. Recall that this is defined to be the minimal number
of $L^2$-balls of radius~$t$ needed to cover~$\calE$.  

\begin{lemma}
\label{lem:suplem} Let $\calE\subset \R_+^N$, $B=\sup_{x\in\calE} |x|$, 
and $\xi_j$ be selector variables as above with $\Prob[\xi_j=1]=\tau$, $\Prob[\xi_j=0]=1-\tau$,
and $0<\tau<1$ arbitrary. Let $1\le m\le N$. Then
\[
\Erw \sup_{x\in\calE,|A|=m} \Big[ \sum_{j\in A} \xi_j\, x_j\Big] \les
(\tau m + 1)^{\half} B + 
\int_0^B \sqrt{\log N_2(\calE,t)}\, dt
\]
where $N_2$ refers to the $L^2$ entropy.
\end{lemma}
\begin{proof} 
Let $\calE_k$ be minimal $2^{-k}$-nets for $\calE$ with $2^{-k}\le B$. Let $B=2^{-k_0}$. 
Then every $x\in\calE$ can be written as
\[ x = x_{k_0} + \sum_{k=k_0}^\infty (x_{k+1}-x_{k}) = x_{k_0}+\sum_{k=k_0}^\infty 2^{-k+1} y_k\]
where $x_k\in\calE_k$ for every $k\ge k_0$. We can and do set $x_{k_0}=0$. 
Now, $y_k\in \calF_k$ where $\diam(\calF_k) \le 1$
and $\#(\calF_k)\le \#(\calE_k)\cdot\#(\calE_{k+1})$. Hence
\begin{equation}
\label{eq:F_kcard}
\log\,\#\calF_k \le C\,\log\,\#\calE_{k+1},
\end{equation}
and thus
\begin{equation}
\label{eq:reduc}
\Erw \sup_{x\in\calE,|A|=m} \Big[ \sum_{j\in A} \xi_j\, x_j\Big]
\le \sum_{k\ge k_0} 2^{-k+1}  \Erw \sup_{y\in\calF_k,\,|A|\le m} \;\sum_{i\in A} \xi_i |y_i|.
\end{equation}
Now fix some $k\ge k_0$ and write $\calF$ instead of $\calF_k$. 
Moreover, replacing every vector $y=\{y_j\}_{j=1}^N\in\calF$ with the vector $\{|y_i|\}_{i=1}^N$,
we may assume that $\calF\in\R_+^N$. Note that this changes neither the diameter nor
the cardinality of $\calF$. 
With $0<\rho_1<\rho_2$ to be determined, one has
\[
\sum_{i\in A} \xi_i y_i \le \sum_{y_i\ge \rho_2} y_i + \sum_{i\in A,\,y_i\le \rho_1} y_i
+\sum_{\rho_1<y_i<\rho_2} \xi_i\, y_i 
\le \rho_2^{-1} \sum_{y_i\ge \rho_2} y_i^2 + m\rho_1 + \sum_{\rho_1<y_i<\rho_2} \xi_i\, y_i. 
\]
Let $q=1+\lfloor\log\calF\rfloor$.
Since $|y|\le1$, one concludes that
\bea
\Erw\sup_{y\in\calF,\,|A|\le m}\sum_{i\in A} \xi_i y_i &\le& \rho_2^{-1} + m\rho_1 + \sup_{y\in\calF} \sum_{\rho_1<y_i<\rho_2} \xi_i\, y_i  \nn \\
&\les& \rho_2^{-1}+m\rho_1 + \Erw \left[\sum_{y\in\calF} \Bigl(\sum_{\rho_1<y_i<\rho_2} \xi_i\, y_i\Bigr)^q \right]^{\frac1q} \label{eq:infty_sum} \\
&\les& \rho_2^{-1}+m\rho_1 + \left[\sum_{y\in\calF} \Erw \Bigl(\sum_{\rho_1<y_i<\rho_2} \xi_i\, y_i\Bigr)^q \right]^{\frac1q} \label{eq:holder} \\
&\les& \rho_2^{-1}+m\rho_1 + (\#\calF)^{\frac1q}\sup_{y\in\calF}\left[ \Erw \Bigl(\sum_{\rho_1<y_i<\rho_2} \xi_i\, y_i\Bigr)^q \right]^{\frac1q} \label{eq:card_out} \\
&\les& \rho_2^{-1}+m\rho_1 + \sup_{|y|\le 1}\Bigl\|\sum_{\rho_1<y_i<\rho_2} \xi_i(\omega)\, y_i\Bigr\|_{L^q(\omega)} \label{eq:sup_F}.
\end{eqnarray}
Here \eqref{eq:infty_sum} follows from the embedding $\ell^q(\calF)\hookrightarrow \ell^\infty(\calF)$,
\eqref{eq:holder} follows from H\"older's inequality, and to pass from \eqref{eq:card_out} to~\eqref{eq:sup_F} one uses that 
\[ (\#\calF)^{\frac1q} = \exp[(\log\#\calF)/q] \le e \]
by our choice of $q=1+\lfloor\log\calF\rfloor$. To control the last term in~\eqref{eq:sup_F}, we need the following simple estimate, see~Lemma~2 in~\cite{B1}. By the multinomial theorem (for any positive integer~$q$),
\bea
\Erw \Big[\sum_{j=1}^n \xi_j\Big]^q &=& \sum_{q_1+\ldots+q_n=q} \binom{q}{q_1,\ldots,q_n} \Erw \xi_1^{q_1}\cdot\ldots \cdot \Erw\xi_n^{q_n} \nn \\
&=&  \sum_{\ell=1}^q \;\sum_{1\le i_1<i_2<\ldots < i_\ell\le n} \;
\sum_{\substack{q_{i_1}+\ldots+q_{i_\ell}=q\\ q_{i_1}\ge 1,\ldots,q_{i_\ell}\ge 1}} \binom{q}{q_{i_1},\ldots,q_{i_\ell}}\tau^\ell \nn \\
&\le& \sum_{\ell=1}^q \frac{n^\ell}{\ell!} \ell^q \tau^\ell \le \sum_{\ell=1}^q \binom{q}{\ell} q^{q-\ell}
\frac{\ell^q}{q!} (n\tau)^\ell \le \sum_{\ell=1}^q \binom{q}{\ell} q^{q-\ell} (e\tau n)^\ell \nn \\
&\le& (q+e\tau n)^q. \label{eq:q_sum}
\end{eqnarray}
It is perhaps more natural (and also more precise) to estimate $q^{\rm th}$ moments by means of
the Bernoulli law
\[
\Erw \Big[\sum_{j=1}^n \xi_j\Big]^q = \sum_{\ell=0}^n \binom{n}{\ell} \ell^q \tau^\ell (1-\tau)^{n-\ell}.
\]
But we have found the approach leading to~\eqref{eq:q_sum} more flexible since it also applies
to non Bernoulli cases.
Continuing with the final term in \eqref{eq:sup_F} one concludes from~\eqref{eq:q_sum} that
\bea
 \sup_{|y|\le 1}\Bigl\|\sum_{\rho_1<y_i<\rho_2} \xi_i(\omega)\, y_i\Bigr\|_{L^q(\omega)} &\le &
2\,\sum_{\rho_2^{-2}<2^j<\rho_1^{-2}} 2^{-\frac{j}{2}} \Bigl\| \sum_{i=1}^{2^j} \xi_i(\omega) \Bigr\|_{L^q(\omega)} 
\label{eq:dyadic} \\
&\le& 2\,\sum_{\rho_2^{-2}<2^j<\rho_1^{-2}} 2^{-\frac{j}{2}} (q+e\tau 2^j) \les q \rho_2 + 
\tau \rho_1^{-1}.
\label{eq:dyad_concl}
\end{eqnarray}
Inserting this bound into~\eqref{eq:sup_F} and setting $\rho_1=\sqrt{\tau/m}$ and $\rho_2=q^{-\half}$
 yields
\[
\Erw\sup_{y\in\calF,\,|A|\le m}\sum_{i\in A} \xi_i y_i \les \sqrt{m\tau} + \sqrt{q} 
\les \sqrt{m\tau} + 1 + \sqrt{\log\#\calF}. 
\]
The lemma now follows in view of \eqref{eq:F_kcard} and~\eqref{eq:reduc}.
\end{proof}

\subsection{Entropy bounds}

As in~\cite{B1} we will need bounds on certain covering numbers, also called entropies.
We recall those bounds starting with the so called ``dual Sudakov inequality''
 for the reader's convenience. More on this can be found in Pisier~\cite{P} and 
Bourgain, Lindenstrauss, Milman~\cite{BLM}, Section~4. Consider $\R^n$ with two norms,
the Euclidean norm $|\cdot|$ and some other (semi)norm~$\|\cdot\|$. 
We set $X=(\R^n,\|\cdot\|)$ and denote the unit ball in this space by~$B_X$,
whereas the Euclidean unit ball will be~$B^n$. As usual, for any set $U\subset\R^n$ 
and $t>0$ one sets
\begin{equation}
\label{eq:Edef}
E(U,B_X,t) := \inf \Big\{ N\ge1\:|\: \exists\;x_j\in\R^n,\,1\le j\le N,\,U\subset \bigcup_{j=1}^N 
(x_j+tB_X) \Big\}. 
\end{equation}
There are two closely related quantities, namely
\bea
\tilde E(U,B_X,t) &:=& \inf \Big\{ N\ge1\:|\: \exists\;x_j\in U,\,1\le j\le N,\,U\subset \bigcup_{j=1}^N (x_j+tB_X) \Big\} \nn \\
D(U,B_X,t) &:=& \sup \Bigl\{ M\ge1\:|\: \exists\; y_j\in U,\, 1\le j\le M,\,\|y_j-y_k\|\ge t,\,j\ne k
\Bigr\}. \label{eq:Ddef}
\end{eqnarray}
There are the following comparisons between these quantities:
\begin{equation}
\label{eq:comp}
D(U,B_X,t)\ge \tilde E(U,B_X,t) \ge E(U,B_X,t) \ge D(U,B_X,2t) \text{\ \ and\ \ }E(U,B_X,t)\ge \tilde
E(U,B_X,2t).
\end{equation}
The final inequality holds because every covering of $U$ by arbitrary $t$-balls gives rise to
a covering by $2t$-balls with centers in~$U$. To see that $E(U,B_X,t) \ge D(U,B_X,2t)$, let 
$\{y_j\}_{j=1}^M\subset U$ be $2t$-separated and $U\subset \bigcup_{i=1}^N (x_i+tB_X)$. Then
every $y_j\in x_{i}+tB_X$ for some $i=i(j)$. Moreover, $j\ne k \Impl i(j)\ne i(k)$. 
Hence $N\ge M$. 

\noindent The ``dual Sudakov inequality'' Lemma~\ref{lem:dual_sud} 
bounds $E(B^n,B_X,t)$ in terms of the Levy mean
\begin{equation}
\label{eq:levy} M_X := \int_{S^{n-1}} \|x\| \,d\sigma(x), 
\end{equation}
where $\sigma$ is the normalized measure on $S^{n-1}$. Alternatively, one has
\bea
M_X &=& \alpha_n (2\pi)^{-\frac{n}{2}}
\int_{\R^n} e^{-\frac{|x|^2}{2}} \;\|x\| \,dx \label{eq:gauss1} \\
M_X &=& \alpha_n \int_{\Omega} \Big\| \sum_{i=1}^n g_i(\omega) \vec{e_i} \Big\| \,d\Prob(\omega), 
\label{eq:gauss2} 
\end{eqnarray}
where
\[ 
\alpha_n = \frac{\Gamma\Big(\frac{n}{2}\Big)}{\Gamma\Big(\frac{n+1}{2}\Big)\sqrt{2}} \asymp n^{-\half}
\]
and $g_i$ are i.i.d.~standard normal variables, and $\vec{e_i}$ is an ONS. The probabilistic form~\eqref{eq:gauss2} is of course just a restatement of~\eqref{eq:gauss1}, whereas the latter can be obtained from the 
definition~\eqref{eq:levy} by means of polar coordinates. 
The following lemma is due to Pajor and Tomczak-Jaegerman~\cite{PT} but the proof given below
 is due to Pajor and Talagrand, see~\cite{BLM}. 

\begin{lemma}
\label{lem:dual_sud}
For any $t>0$
\begin{equation}
\label{eq:dual_sud} \log E(B^n,B_X,t) \le Cn \Big(\frac{M_X}{t}\Big)^2,
\end{equation}
where $C$ is an absolute constant.
\end{lemma}
\begin{proof}
Let $\{x_i\}_{i=1}^N \subset B^n$, $\|x_i-x_j\|\ge t$ for $i\ne j$ and $N$ maximal. 
Then $E(B^n,B_X,t)\le N$. Let $\mu(dx)=(2\pi)^{-\frac{n}{2}} e^{-\frac{|x|^2}{2}} \; dx$.
Then by definition~\eqref{eq:levy},
\begin{equation}
\label{eq:median}
 \mu(\|x\|>2M_X\alpha_n^{-1})<\half \Impl \mu(\|x\|\le 2M_X\alpha_n^{-1})>\half.
\end{equation}
Moreover, $\{x_i+\half tB_X\}_{i=1}^N$ and therefore also 
$\{ y_i+2M_X \alpha_n^{-1} B_X\}_{i=1}^N$ have mutually disjoint interiors,
where we have set $y_i=4M_X (t\alpha_n)^{-1} x_i$.  Now, by symmetry of $B_X$ and convexity of~$e^{-u}$,
\bea
\mu(y_i+2M_X \alpha_n^{-1} B_X) &=& (2\pi)^{-\frac{n}{2}} \int_{2M_X \alpha_n^{-1} B_X} 
e^{-|y-y_i|^2/2} \, dy \nn \\
&=& (2\pi)^{-\frac{n}{2}} \int_{2M_X \alpha_n^{-1} B_X} 
\half\left[e^{-|y-y_i|^2/2} +e^{-|y+y_i|^2/2}\right] \, dy \nn \\
&\ge& (2\pi)^{-\frac{n}{2}} \int_{2M_X \alpha_n^{-1} B_X} 
e^{-(|y-y_i|^2+|y+y_i|^2)/4} \, dy \nn \\
&=& (2\pi)^{-\frac{n}{2}} \int_{2M_X \alpha_n^{-1} B_X} e^{-(|y|^2+|y_i|^2)/2} \, dy \ge \half e^{-|y_i|^2/2}, \nn
\end{eqnarray}
where the last step follows from~\eqref{eq:median}. Since $|y_i|\le 4M_X(t\alpha_n)^{-1}$,
\[ \mu(y_i+2M_X \alpha_n^{-1} B_X) \ge \half \exp\Big(-\half (4M_X)^2(t\alpha_n)^{-2}\Big).\]
Hence
\[ 1\ge \sum_{i=1}^N \mu(y_i+2M_X \alpha_n^{-1} B_X) \ge \half N \exp\Big(-(4M_X)^2(t\alpha_n)^{-2}\Big),\]
and the lemma follows since $\alpha_n\asymp n^{-\half}$.
\end{proof}

Observe that \eqref{eq:dual_sud} is a poor bound as $t\to0$. Indeed, rather than
the~$\exp(t^{-2})$ behavior exhibited by~\eqref{eq:dual_sud} the true asymptotics is~$t^{-n}$  
as~$t\to0$. The point of Lemma~\ref{lem:dual_sud} is to relate the size of~$t$ to both
$M_X$ and~$n$. This is best illustrated by some standard examples.
\begin{itemize}
\item
Firstly, take $X=\ell^1_n$. In that case, 
\[
\alpha_n^{-1} M_X = (2\pi)^{-\frac{n}{2}} \int_{\R^n} \sum_{i=1}^n |x_i| e^{-\frac{|x|^2}{2}}\,dx 
= \frac{n}{\sqrt{2\pi}} \int_{-\infty}^\infty |x_1| e^{-\frac{x_1^2}{2}} \, dx_1 = \frac{2n}{\sqrt{2\pi}}. \nn
\]
Therefore, $M_X\asymp \sqrt{n}$. By~\eqref{eq:dual_sud}, 
\[  \sup_n E(B^n,B_{\ell^1_n},n) \le C.\] 
This bound is somewhat wasteful. Indeed, since $\sqrt{n}B_{\ell^1_n} \supset B^n$, one actually has
\[  \sup_n E(B^n,B_{\ell^1_n},\sqrt{n}) \le C.\] 
The reason for this ``overshoot'' is that the major contribution to $M_X$ comes from the
corners of $B_{\ell^1_n}$. On the other hand, these corners do not determine the smallest $r$ for 
which~$rB_X\supset B^n$.
\item Secondly, consider $X=\ell^\infty_n$. Using~\eqref{eq:gauss2}, 
\[ \alpha_n^{-1}M_X = \Erw \sup_{1\le i\le n} |g_i| \asymp \sqrt{\log n},\]
where the latter bound is a rather obvious and well-known fact. Hence
\[ M_X \asymp \sqrt{\frac{\log n}{n}}\]
which implies via~\eqref{eq:dual_sud} that
\[ \sup_n E(B^n,B_{\ell^\infty_n},\sqrt{\log n}) \le C.\] 
This is the correct behavior up to the $\log n$-factor since $B^n \subset B_{\ell^\infty_n}$. 
In contrast to the previous case, the bulk of the contribution to $M_X$ comes from that part
of $B_{\ell^\infty_n}$ that is also the most relevant for the covering of the Euclidean ball.
\item Finally, and most relevantly for our purposes, identify $\R^n$ with the space of 
trigonometric polynomials with  real coefficients of degree $n$, i.e.,
\begin{equation}
\label{eq:iden}
\R^n \simeq \Bigl\{ \sum_{j=1}^n a_j e(j\theta) \:|\: a_j\in \R \Bigr\}. 
\end{equation}
Furthermore, define $\|\cdot\|= \|\cdot\|_{L^q(\tor)}$ where $q\ge2$ is fixed. 
Then
\bea
 M_X &=& \alpha_n \int_{\Omega} \Big\| \sum_{j=1}^n g_j(\omega) e(j\theta) \Bigr\|_{L^q(\tor)} \, d\Prob(\omega) \nn \\
&=& \alpha_n \Erw \int_{\Omega} \Big\| \sum_{j=1}^n \pm g_j(\omega) e(j\theta) \Bigr\|_{L^q(\tor)} \, d\Prob(\omega) \label{eq:rade} \\
&\le& C \alpha_n \sqrt{q} \int_{\Omega} \Big(\sum_{j=1}^n g_j^2(\omega) \Bigr)^{\half} \, d\Prob(\omega) 
\label{eq:khinchin} \\
&\le& C\alpha_n \sqrt{q} \left(\int_{\Omega} \sum_{j=1}^n g_j^2(\omega) \, d\Prob(\omega) \right)^{\half} 
= C\alpha_n \sqrt{q} \sqrt{n} \le C\sqrt{q}. \nn
\end{eqnarray}
In \eqref{eq:rade} the expectation $\Erw$ refers to the random and symmetric choice of signs $\pm$,
whereas the $\sqrt{q}$-factor in~\eqref{eq:khinchin} is due to the fact that the constant in 
Khinchin's inequality grows like~$\sqrt{q}$. Hence
\begin{equation}
\label{eq:Lqentropy} \log E(B^n, B_X,t) \le C\,qn t^{-2} 
\end{equation}
in this case. 
\end{itemize} 

\noindent
The proof of Theorem~\ref{thm:Lp} requires estimating $N_q(\calP_A,t):= E(\calP_A,B_{L^q(\tor)},t)$. Here 
\[ 
\calP_A := \Big\{ \sum_{n\in A} a_n e(n\theta) \:\Big|\: |a|=|a|_{\ell^2_N} \le 1 \Big\}
\]
where $ A\subset [1,N]$. 
Invoking~\eqref{eq:Lqentropy} leads to 
\begin{equation}
\label{eq:easy_ent}
 \log N_q(\calP_A,t) \le Cq|A|\,t^{-2}.
\end{equation}
This bound is basically optimal when $t\sim 1$, but
it can be improved for very small and very large~$t$. 

\begin{cor}
\label{cor:small_t}
For $q\ge2$ and any $ A\subset [1,N]$
\begin{equation}
\log N_q(\calP_A,t) \le C\,q|A|\Big[1+ \log\frac{1}{t}\Big] \text{\ \ if\ \ }0<t\le\half. 
\label{eq:klein_t}
\end{equation}
\end{cor}
\begin{proof}
Let $m=|A|$. Thus $1\le m\le N$. 
Notice firstly that
\begin{equation}
\label{eq:same_cover} 
\log N_q\Big(\Big\{\sum_{n\in A} a_n e(n\theta)\:\Big|\: |a|\le 1\Big\},t\Big) 
\le C\,m\log \frac{1}{t} + \log N_q\Big(\Big\{\sum_{n\in A} a_n e(n\theta)\:\Big|\: |a|\le 1\Big\},1\Big).
\end{equation}
This follows from the fact that for any norm $\|\cdot\|$ in $\R^m$ with unit-balls $B_X$ one has
\begin{equation}
\label{eq:D_scale}
 D(B_X,B_X,t) \le (4/t)^m \text{\ \ for all\ \ }0<t<1
\end{equation}
by scaling and volume counting, see~\eqref{eq:Ddef} for the definition of $D(B_X,B_X,t)$. 
Indeed, suppose $M=D(B_X,B_X,t)$. Then there are $M$ disjoint balls 
$\{x_j+\half tB_X\}_{j=1}^M$ with centers $x_j\in B_X$. Since $x_j+\half tB_X \subset 2B_X$ if $t<1$, 
it follows that
\[ \sum_{j=1}^M |\half t B_X| \le |2B_X| \Impl M (t/2)^m \le 2^m, \]
as claimed. Here $|\cdot|$ stands for Lebesgue measure. Thus~\eqref{eq:D_scale} holds, and
therefore also~\eqref{eq:same_cover} in view of~\eqref{eq:comp}. Hence
\bea
\log N_q(\calP_A,t) &\le&   C\,m\log\frac{1}{t} + 
\log N_q\Big(\Big\{\sum_{n\in A} a_n e(n\theta)\:\Big|\: |a|\le 1\Big\},1\Big) \nn \\
&\le&  C\,m\log\frac{1}{t} + Cqm, \nn
\end{eqnarray}
where the final term follows from~\eqref{eq:Lqentropy}. 
\end{proof}

\noindent We now turn to large $t$. 
The following corollary slightly improves on the rate of decay.

\begin{cor}
\label{cor:entropy2}
Let $q\ge2$ and $A\subset [1,N]$. 
With $\calP_A$ as above one has
\begin{equation}
\log N_q(\calP_A,t) \le Cq\,|A|\, t^{-\nu} \text{\ \ if\ \ }t>\half \label{eq:gross_t} 
\end{equation}
where $\nu=\nu(q)>2$.
\end{cor}
\begin{proof}
Recall that $N_q(\calP_A,t)=E(\calP_A,B_{L^q},t)$. Using~\eqref{eq:comp}, one obtains 
from~\eqref{eq:easy_ent} that also
\begin{equation}
\label{eq:centersinE} 
\log \tilde E(\calP_A,B_{L^q},t) \le  Cq\,|A|\,t^{-2}.
\end{equation}
Let $q<r$, $\frac{1}{q}=\frac{1-\theta}{2}+\frac{\theta}{r}$. Since for any $f,g\in \calP_A$ 
\[ \|f-g\|_q \le \|f-g\|_2^{1-\theta}\|f-g\|_r^\theta \le 2\|f-g\|_r^\theta,\]
one concludes from \eqref{eq:centersinE} that
\[ 
\log \tilde E(\calP_A,B_{L^q},t) \le \log \tilde E(\calP_A,B_{L^r},(t/2)^{1/\theta}\big) \le
Cq\,|A|\,t^{-2/\theta}.
\]
Applying~\eqref{eq:comp} again yields~\eqref{eq:gross_t}.
\end{proof}

\subsection{Decoupling lemma}

Lastly, we require a version of Bourgain's decoupling technique, cf.~Lemma~4 in~\cite{B1}.
In contrast to his case we only need to decouple into two sets rather than three.

\begin{lemma}
\label{lem:decoup}
Let real-valued functions $h_\alpha(u)$ on $\R$ be given for $\alpha=1,2,3$ that satisfy 
\[ |h_\alpha(u)|\le C(1+|u|)^{p_\alpha},\quad |h_\alpha(u)-h_\alpha(v)|\le C(1+|u|+|v|)^{p_\alpha-\delta}|u-v|\]
for all $u,v\in\R$ and some fixed choice of $p_\alpha>0$, $\delta>0$. 
Let $x,y,z\in\ell^2_N$ be sequences 
so that $|x|,|y|,|z|\le 1$ and suppose $\zeta_j=\zeta_j(t)$ are i.i.d.~random variables with
$\Prob(\zeta_j=1)=\Prob(\zeta_j=0)=\half$. We assume that $\Prob(dt)=dt$ on $[0,1]$, say. Set 
$ R^1_t = \{1\le j\le N\:|\: \zeta_j(t)=1\}, \quad R^2_t = \{1\le j\le N\:|\: \zeta_j(t)=0\}.$
Then
\bea
&& \left| \int h_1\Big( \sum_{i\in R^1_t} x_i\Big) h_2\Big( \sum_{i\in R^2_t} y_i\Big)
h_3\Big( \sum_{i\in R^2_t} z_i\Big) \,dt - h_1\Big( \half\sum_{i} x_i\Big) h_2\Big( \half\sum_{i} y_i\Big)
h_3\Big( \half\sum_{i} z_i\Big) \right| \nn \\
&& \le C\,\Bigl(1+ \Big|\sum_{i} x_i\Big| + \Big|\sum_{i} y_i\Big|+  \Big|\sum_{i} z_i\Big| \Bigr)^{p-\delta}
\end{eqnarray}
where $p=p_1+p_2+p_3$ and $C$ is some absolute constant. 
\end{lemma}  
\begin{proof}
By assumption,
\bea
\Big|h_\alpha\Big( \sum_{i\in R^1_t} x_i\Big) - h_\alpha\Big(\half \sum_{i=1}^N x_i\Big)\Big|
&\le& C\Big(1+\big|\sum_{i=1}^N x_i\big|+\big|\sum_{i=1}^N (\zeta_i-\half)x_i\big|\Big)^{p_\alpha-\delta}
\Big| \sum_{i=1}^N (\zeta_i-\half)x_i \Big|^{\delta} \nn \\
&\le& C\Big(1+\big|\sum_{i=1}^N x_i\big|\Big)^{p_\alpha-\delta}\Big(1+\big|\sum_{i=1}^N (\zeta_i-\half)x_i\big|\Big)^{p} \nn \\
\Big|h_\alpha\Big( \sum_{i\in R^1_t} x_i\Big)\Big|+ \Big|h_\alpha\Big(\half \sum_{i=1}^N x_i\Big)\Big|
&\le& C\Big(1+\big|\sum_{i=1}^N x_i\big|\Big)^{p_\alpha}\Big(1+\big|\sum_{i=1}^N (\zeta_i-\half)x_i\big|\Big)^{p} \nn 
\end{eqnarray}
for $\alpha=1,2,3$. Hence
\bea
&& \left| \int h_1\Big( \sum_{i\in R^1_t} x_i\Big) h_2\Big( \sum_{i\in R^2_t} y_i\Big)
h_3\Big( \sum_{i\in R^2_t} z_i\Big) \,dt - h_1\Big( \half\sum_{i=1}^N x_i\Big) h_2\Big( \half\sum_{i=1}^N y_i\Big)
h_3\Big( \half\sum_{i=1}^N z_i\Big) \right| \nn \\
&& \le C\Big(1+ \Big|\sum_{i=1}^N x_i\Big| + \Big|\sum_{i=1}^N y_i\Big|+  \Big|\sum_{i=1}^N z_i\Big| \Bigr)^{p-\delta}
\nn \\
&& \mbox{\hspace{1in}} \cdot \int \left(1+\Big| \sum_{i=1}^N (\zeta_i-\half)x_i \Big|+\Big| \sum_{i=1}^N (\zeta_i-\half)y_i \Big|
+ \Big| \sum_{i=1}^N (\zeta_i-\half)z_i \Big|\right)^{p} \, dt.
\end{eqnarray}
The lemma now follows from Khinchin's inequality. Indeed,
\[ \int \Big| \sum_{i=1}^N (\zeta_i-\half)x_i \Big|^p \,dt \le C_p\,|x|^p\le C_p,\]
by assumption.
\end{proof}

\subsection{The proof of Theorem~\ref{thm:Lp} for $p=3$}

We now start the proof of Theorem~\ref{thm:Lp} for $p=3$. 
In fact, we state a somewhat more precise form of this theorem for $p=3$.

\begin{theorem}
\label{thm:L3} 
Let $0<\delta<1$ be fixed. For every positive integer $N$ we let $\xi_j=\xi_j(\omega)$ be i.i.d.~variables with $\Prob[\xi_j=1]=\tau$, $\Prob[\xi_j=0]=1-\tau$ where $\tau=N^{-\delta}$. Define a random subset
\[ S(\omega)=\{j\in[1,N]\:|\: \xi_j(\omega)=1\}.\]
Then for every $\gamma>0$ there is a constant $C_\gamma$ so that 
\begin{equation}
\label{eq:prob3} 
\sup_{N\ge1}\;\Prob\Bigl[\sup_{|a_n|\le1} \Big\|\sum_{n\in S(\omega)} a_n e(n\theta) \Big\|_{L^3(\tor)} \ge C_\gamma\, \Big\|\sum_{n\in S(\omega)}  e(n\theta) \Big\|_{L^3(\tor)} \Bigr]  \le \gamma.
\end{equation}
\end{theorem}
\begin{proof}
Firstly, note that for fixed $0<\delta<1$ and large $N$ Lemma~\ref{lem:kinchin} implies that
\[ \Prob\Big[\summe \xi_n \ge 2\tau N\Big] \les \exp(-c\tau N). \]
Let $\Erw'$ denote the restricted expectation
\[ 
\Erw' \sup_{|a_n|\le1} \Big\|\summe\xi_n a_n e(n\theta) \Big\|_{L^3(\tor)} := \Erw \chi_{[\sum \xi_n \le 2\tau N]}\;\sup_{|a_n|\le1} \Big\|\summe\xi_n a_n e(n\theta) \Big\|_{L^3(\tor)}. 
\]
Then
\bea
\Erw \sup_{|a_n|\le1} \Big\|\summe\xi_n a_n e(n\theta) \Big\|_{L^3(\tor)} &\le& 
N\exp(-c\tau N)+\Erw' \sup_{|a_n|\le1} \Big\|\summe\xi_n a_n e(n\theta) \Big\|_{L^3(\tor)} \nn \\
&\le& O(1) + \Erw' \sup_{|a_n|\le1} \Big\|\summe\xi_n a_n e(n\theta) \Big\|_{L^3(\tor)}. \nn
\end{eqnarray} 
From now on, we set $m=2\tau N$, and we will mostly work with $\Erw'$ instead of $\Erw$.
Next, fix some $\{a_n\}_{n=1}^N$ with $|a_n|\le 1$. Then, rescaling Lemma~\ref{lem:decoup} (with $h_1(x)=h_2(x)=x$ and $h_3(x)=|x|$)  one obtains that
\bea
\frac18 \int_0^1 \left|\summe a_n \xi_n e(n\theta) \right|^3 \,d\theta &=& 
\int \int_0^1 \sum_{n\in R^1_t} a_n\xi_n e(n\theta) \sum_{k\in R^2_t} \bar a_k\xi_k e(-k\theta) 
\Big|\sum_{\ell\in R^2_t} a_\ell\xi_\ell e(\ell\theta) \Big|\, d\theta\,dt \nn \\
&& + O\left( m^{\frac32} \int_0^1 \Big(1+\Big|\summe \frac{a_n}{\sqrt{m}}\xi_n e(n\theta)\Bigr|^2 \Big)\, d\theta \right). \label{eq:step1}
\end{eqnarray}
The $O$-term in \eqref{eq:step1} is $O(m^{\frac32})$ by construction.  Let $\{\xi_n(\omega_1)\}_{n=1}^N$ and $\{\xi_n(\omega_2)\}_{n=1}^N$ 
denote two independent copies of~$\{\xi_n(\omega)\}_{n=1}^N$. Recall that $R_t^1$ and $R_t^2$ are
disjoint for every~$t$. Therefore, {\em for fixed}~$t$ 
\bea
&& \Erw_{\omega} \sup_{|a_n|\le 1} \left| \int_0^1  \sum_{n\in R^1_t} a_n\xi_n(\omega) e(n\theta) \sum_{k\in R^2_t} \bar a_k\xi_k(\omega) e(-k\theta) 
\Big|\sum_{\ell\in R^2_t} a_\ell\xi_\ell(\omega) e(\ell\theta) \Big|\, d\theta \; \right| \nn \\
&& = \Erw_{\omega_1,\omega_2} \sup_{|a_n|\le 1} \left| \int_0^1  \sum_{n\in R^1_t} a_n\xi_n(\omega_1) e(n\theta) \sum_{k\in R^2_t} \bar a_k\xi_k(\omega_2) e(-k\theta) 
\Big|\sum_{\ell\in R^2_t} a_\ell\xi_\ell(\omega_2) e(\ell\theta) \Big|\, d\theta \; \right|. \label{eq:exp_split}
\end{eqnarray}
This leads to 
\bea
&& \Erw_\omega \sup_{|a_n|\le 1} \int_0^1 \Big| \summe a_n\xi_n(\omega) e(n\theta) \Big|^3 \, d\theta \nn \\
&\les& m^{\frac32} 
+ \int \Erw_{\omega_1,\omega_2} \sup_{|a_n|\le 1} \left| \int_0^1  \sum_{n\in R^1_t} a_n\xi_n(\omega_1) e(n\theta) \sum_{k\in R^2_t} \bar a_k\xi_k(\omega_2) e(-k\theta) 
\Big|\sum_{\ell\in R^2_t} a_\ell\xi_\ell(\omega_2) e(\ell\theta) \Big|\, d\theta \; \right|\,dt \nn \\
&\les& m^{\frac32} 
+ \int \Erw_{\omega_1}'\Erw_{\omega_2}' \sup_{\substack{|a_n|\le 1\\|b_n|\le1}} \left| \int_0^1  \sum_{n\in R^1_t} a_n\xi_n(\omega_1) e(n\theta) \sum_{k\in R^2_t} \bar b_k\xi_k(\omega_2) e(-k\theta) 
\Big|\sum_{\ell\in R^2_t} b_\ell\xi_\ell(\omega_2) e(\ell\theta) \Big|\, d\theta \; \right|\,dt \nn \\
&\les& m^{\frac32} + \int \Erw_{\omega_1}'\Erw_{\omega_2}' \sup_{\substack{|a_n|\le 1\\|b_n|\le1}}
\left|\int_0^1  \sum_{n=1}^N a_n\xi_n(\omega_1) e(n\theta) \sum_{k=1}^N \bar b_k\xi_k(\omega_2) e(-k\theta) 
\Big|\sum_{\ell=1}^N b_\ell\xi_\ell(\omega_2) e(\ell\theta) \Big|\, d\theta \; \right|\,dt \nn \\
&& \quad \les m^{\frac32} + \Erw'_{\omega_2} \Erw_{\omega_1}  \sup_{x\in \calE(\omega_2)} \sup_{|A|= m} 
\sum_{n\in A} \xi_n(\omega_1) x_n. \label{eq:step2}
\end{eqnarray}
Here
\[ 
\calE(\omega_2) := \Bigl\{ \Bigl(\Bigl|\Bigl\la e(n\cdot), 
\sum_{k=1}^N \bar b_k\xi_k(\omega_2) e(-k\cdot)\Big|\sum_{\ell=1}^N b_\ell\xi_\ell(\omega_2) e(\ell\cdot) \Big| \Bigr \ra \Bigr| \Bigr)_{n=1}^N \:\Big|\: \sup_{1\le n\le N}|b_n|\le 1\Bigr\} \subset \R^N_+.
\]
In the calculation leading up to~\eqref{eq:step2} we firstly used~\eqref{eq:exp_split}, 
secondly the obvious fact that the supremum only increases if we introduce~$\{b_n\}_{n=1}^N$ 
in addition to~$\{a_n\}_{n=1}^N$, thirdly that one can remove the restrictions to the sets~$R^1_t$ 
and~$R^2_t$ because they can be absorbed into the choice of the sequences~$a_n, b_n$, and lastly 
that $\sum_n \xi_n \le m$ which allows us to introduce~$A\subset [1,N]$, $|A|= m$. 
If $x\in \calE(\omega_2)$, then
\begin{equation}
\label{eq:L4maj} 
|x|_{\ell^2_N}^2 \le \sup_{|a_k|\le1} \Big\|\sum_{k} a_k\,\xi_k (\omega_2) e(k\cdot)\Bigr\|_4^4
\le \Big\|\sum_{k} \xi_k (\omega_2) e(k\cdot)\Bigr\|_4^4 =: B_4^2(\omega_2)
\end{equation}
by the $L^4$ majorant property. By Lemma~\ref{lem:norm1},
\begin{equation}
\label{eq:erwB4}
 \Erw B_4 \le \Bigl(\Erw I_{4,N}\Bigr)^{\half} \les \tau^2 N^{\frac32} + \tau N.
\end{equation}
We now apply Lemma~\ref{lem:suplem} to \eqref{eq:step2}. This yields
\bea
&& \Erw_\omega \sup_{|a_n|\le 1} \int_0^1 \Big| \summe a_n\xi_n(\omega) e(n\theta) \Big|^3 \, d\theta  \les m^{\frac32}\! + \!\Erw'_{\omega_2} \left[ (\sqrt{\tau m}+1)\, B_4(\omega_2)\! + \!\! \int_0^\infty \sqrt{\log N_2(\calE(\omega_2),t)}\,dt\right] \nn \\
&& \qquad \les (\tau N)^{\frac32} +(1+\tau N^{\half}) (\tau^2 N^{\frac32} + \tau N) + \Erw'_{\omega_2} \int_0^\infty \sqrt{\log N_2(\calE(\omega_2),t)}\,dt. \label{eq:step3}
\end{eqnarray}
It remains to deal with the entropy integral in~\eqref{eq:step3}. To this end, observe that the
distance between any two elements in $\calE(\omega_2)$ is of the form
\bea
 \|g|g|-h|h|\|_2 &\les& \|g-h\|_\infty (\|g\|_{2} + \|h\|_{2}) \nn \\
&\les& N^\eps \|g-h\|_q (\|g\|_{2} + \|h\|_{2}) \les N^\eps \sqrt{m} \|g-h\|_q,\nn
\end{eqnarray}
where we chose $q$ very large depending on $\eps$ 
(the factor $N^\eps$ comes from Bernstein's inequality). 
Here $g,h\in \sqrt{m}\calP_A$ where $A=A(\omega_2)=\{n\in[1,N]\:|\: \xi_n(\omega_2)=1\}$ and 
\begin{equation}
\label{eq:PAdef} \calP_A = \left\{ \sum_{n\in A} a_n e(n\cdot)\:\Big|\: |a|_{\ell^2_N}\le 1\right\}.
\end{equation}
Actually, our coefficients are in the unit-ball of $\ell^\infty_n$, but
we have embedded this into $\ell^2_m$ in the obvious way, which leads to the $\sqrt{m}$-factor in front of $\calP_A$ (at this point recall that we are working with $\Erw_{\omega_2}'$). One concludes that, for $\eps>0$ small and $q<\infty$ large depending on $\eps$, 
\bea
\log N_2(\calE(\omega_2),t) &\le& \log N_q(\calP_A, N^{-\eps} m^{-1}t) \nn \\
&\le& Cq\,m \left\{ \begin{array}{ll} 1 + \log \frac{1}{t} & \quad 0<t<m N^\eps \\
                                    (m^{-1} N^{-\eps} t)^{-\nu} & \quad t> N^\eps m
                 \end{array}
         \right.
\end{eqnarray}
where $\nu>2$, see Corollary~\ref{cor:small_t} and Corollary~\ref{cor:entropy2}.  
It follows that the last term in~\eqref{eq:step3} is at most 
\[
\Erw_{\omega_2} \int_0^\infty \sqrt{\log N_2(\calE(\omega_2),t)}\,dt 
\les N^\eps m^{\frac32}.
\]
Plugging this into~\eqref{eq:step3} yields
\bea
 \Erw_\omega \sup_{|a_n|\le 1} \int_0^1 \Big| \summe a_n\xi_n(\omega) e(n\theta) \Big|^3 \, d\theta  &\les&   (\tau N)^{\frac32} +(1+\tau N^{\half}) (\tau^2 N^{\frac32} + \tau N) + N^\eps (\tau N)^{\frac32} \nn \\
&\les& \tau^3 N^2 + N^\eps (\tau N)^{\frac32}. \label{eq:upp_eps}
\end{eqnarray}
Now suppose $\delta< \frac13$. Then $\tau^3 N^2>N^\eps (\tau N)^{\frac32}$ provided $\eps>0$
is small and fixed, and provided $N$ is large. Hence, combining~\eqref{eq:upp_eps}
 with Lemma~\ref{lem:norm1} leads to Theorem~\ref{thm:L3} at least if $\delta<\frac13$.
If one is willing to loose a $N^\eps$-factor, then~\eqref{eq:upp_eps} in combination with Lemma~\ref{lem:norm1} leads to the desired bounds in all cases. 
On the other hand, if $\delta\ge \frac13$ so that typically $\#(S(\omega))\les N^{\frac23}$,
then Bourgain showed that $S(\omega)$ is a $\Lambda_3$ set with large probability.
More precisely, he showed that the constant
\[ K_3(\omega):=\sup_{|a|_{\ell^2_N}\le 1} \Big\|\sum_{n\in S(\omega)} a_n e(n\cdot) \Big\|_3\]
satisfies $\Erw K_3^3 \le C$, see also Theorem~\ref{thm:jean} below. Hence, in our case,
\[ \Erw  \sup_{|a_n|\le 1}  \Big\| \summe a_n\xi_n(\omega) e(n\cdot) \Big\|_3^3 \les (\tau N)^{\frac32}.\]
Clearly,
\[  \Big\| \summe \xi_n(\omega) e(n\cdot) \Big\|_3 \ge \Big\| \summe \xi_n(\omega) e(n\cdot) \Big\|_2 =
\#(S(\omega))^{\half},\]
and we have thus proved \eqref{eq:prob3} for $\delta\ge \frac13$ as well.
\end{proof}

\noindent 
It is perhaps worth pointing out that interpolation of the $L^4$ bound with the $L^2$ bound gives
\[ \tau^{\frac52} N^2 + (\tau N)^{\frac32},\]
so that the estimate we just obtained is better by the initial $\tau^3$-factor (note that this is
due to the $\sqrt{\tau m}$-factor in Lemma~\ref{lem:suplem} as compared to a $\sqrt{\tau N}$-factor). 

\subsection{The case of general $p$}

The strategy is to first generalize the previous argument to all odd integers using
the fact that the majorant property holds for all even integers (for $p=3$ we used this fact with $p=4$). 
Then one runs the same argument again, using now that the (random) majorant property holds for {\it all} 
integers~$p$ and so on. For a given $\eps>0$ this yields that there is a set of $p$ that is 
$\eps$-dense in $[2,\infty)$ and for which the majorant property holds. This is enough by interpolation, 
since we are allowing a loss of~$N^\eps$ in~\eqref{eq:prob}. 
Unfortunately, there are certain technical complications in carrying out this program having to
do with the size of~$\delta$. In this section we deal with~$\delta\le\half$, and in the following
section we discuss a refinement of the method that allows one to relax this condition in some cases.

\noindent Lemma~\ref{lem:main_step} formalizes the main probabilistic argument from the previous section.
Let $p\ge2$. In this section, we say that the {\it random majorant property} (or RMP in short) 
holds at~$p$ if and only if for every $\eps>0$ there exists a constant $C_\eps$ so that 
\begin{equation}
\label{eq:rmp}
\Erw \sup_{|a_n|\le1} \Big\|\sum_{n=1}^N a_n\xi_n e(n\theta) \Bigr\|_p^p \le C_\eps\,N^\eps 
\Erw  \Big\|\sum_{n=1}^N \xi_n e(n\theta) \Bigr\|_p^p
\end{equation}
for all $N\ge1$. Note that (the proof of) Theorem~\ref{thm:L3} establishes 
that the random majorant property holds at~$p=3$. Moreover, if~\eqref{eq:rmp} holds
for some~$p$, then~\eqref{eq:prob} also holds for that value of~$p$, see Lemma~\ref{lem:norm1}.

\begin{lemma}
\label{lem:main_step}
Let $2\le p\le 3$. Suppose the random majorant property~\eqref{eq:rmp} holds at~$2(p-1)$. Then
it also holds at~$p$. Furthermore, 
suppose the RMP holds at $p-1$, $2(p-1)$ and $2(p-2)$. If $4\ge p\ge3$, then
it also holds at~$p$. If $p>4$ and $\delta\le \half$ (i.e., $\tau=N^{-\delta}\ge N^{-\half}$), 
then it also holds at~$p$.
\end{lemma}
\begin{proof}
Assume first that $p\ge3$. Instead of~\eqref{eq:step1}, Lemma~\ref{lem:decoup} implies in this case that
\bea
2^{-p} \int_0^1 \left|\summe a_n \xi_n e(n\theta) \right|^p \,d\theta &=& 
\int \int_0^1 \sum_{n\in R^1_t} a_n\xi_n e(n\theta) \sum_{k\in R^2_t} \bar a_k\xi_k e(-k\theta) 
\Big|\sum_{\ell\in R^2_t} a_\ell\xi_\ell e(\ell\theta) \Big|^{p-2}\, d\theta\,dt \nn \\
&& + O\left( m^{\frac{p}{2}} \int_0^1 \Big(1+\Big|\summe \frac{a_n}{\sqrt{m}}\xi_n e(n\theta)\Bigr|^{p-1} \Big)\, d\theta \right). \label{eq:step1_p}
\end{eqnarray}
To bound the $O$-term in \eqref{eq:step1_p} note that by the RMP for $p-1\ge 2$, 
\begin{equation}
\Erw \sup_{|a_n|\le1}\int_0^1 \Big|\summe a_n\xi_n e(n\theta)\Bigr|^{p-1} \,d\theta \le
C_\eps\,N^\eps\Erw \int_0^1 \Big|\summe \xi_n e(n\theta)\Bigr|^{p-1} \,d\theta 
= C_\eps\,N^\eps \Erw I_{p-1,N}. \label{eq:I_{p-1}1}
\end{equation}
A calculation analogous to that leading up to~\eqref{eq:step2} therefore yields
\begin{equation}
\Erw_\omega \sup_{|a_n|\le 1} \int_0^1 \Big| \summe a_n\xi_n(\omega) e(n\theta) \Big|^p \, d\theta 
\les m^{\frac{p}{2}} + C_\eps\,N^\eps m^{\half}  \Erw I_{p-1,N} 
 + \Erw'_{\omega_2} \Erw'_{\omega_1}  \sup_{x\in \calE(\omega_2)} \sup_{|A|=m} 
\sum_{n\in A} \xi_n(\omega_1) x_n, \label{eq:step2_p}
\end{equation}
where now
\[ 
\calE(\omega_2) = \Bigl\{ \Bigl(\Bigl|\Bigl\la e(n\cdot), 
\sum_{k=1}^N \bar b_k\xi_k(\omega_2) e(-k\cdot)\Big|\sum_{\ell=1}^N b_\ell\xi_\ell(\omega_2) e(\ell\cdot) \Big|^{p-2} \Bigr \ra \Bigr| \Bigr)_{n=1}^N \:\Big|\: \sup_{1\le n\le N}|b_n|\le 1\Bigr\} \subset \R^N_+.
\]
If $x\in \calE(\omega_2)$, then by Plancherel and the RMP at $2(p-1)$, 
\bea
\label{eq:L2kmaj} 
\Erw \sup_{x\in\calE(\omega_2)} |x|_{\ell^2_N}^2 &\le& \Erw_{\omega_2}\sup_{|a_k|\le1} \int_0^1 \Big|\sum_{k} a_k\,\xi_k (\omega_2) e(k\theta)\Bigr|^{2(p-1)}
\,d\theta \\
&\le& C_\eps\,N^\eps \Erw_{\omega_2}\; \int_0^1 \Big|\sum_{k} \xi_k (\omega_2) e(k\theta)\Bigr|^{2(p-1)}\,d\theta 
\le C_\eps\, N^\eps \Erw I_{{2(p-1)},N}. \nn 
\end{eqnarray}
Thus, by \eqref{eq:step2_p} and Lemma~\ref{lem:suplem},
\bea
&&\Erw_\omega \sup_{|a_n|\le 1} \int_0^1 \Big| \summe a_n\xi_n(\omega) e(n\theta) \Big|^p \, d\theta \nn\\
&& \le C_\eps\,N^\eps\Big[m^{\frac{p}{2}}+m^{\frac12} \Erw I_{p-1,N} + (1+\sqrt{m\tau}) \sqrt{\Erw I_{2(p-1),N}} +
\Erw_{\omega_2}' \int_0^\infty \sqrt{\log N_2(\calE(\omega_2),t)}\,dt\Big]. 
\label{eq:step4_p}
\end{eqnarray}
To estimate the entropy term, let $q$ be very large depending on~$\eps$. 
Then the distance between any two elements in $\calE(\omega_2)$ is of the form
\bea
 \|g|g|^{p-2}-h|h|^{p-2}\|_2 &\le& \|g-h\|_\infty \big(\|g\|_{2(p-2)}^{p-2} + \|h\|_{2(p-2)}^{p-2}\big) \nn \\
&\le& C_\eps\,N^{\eps} \|g-h\|_q \big(\|g\|_{2(p-2)}^{p-2} + \|h\|_{2(p-2)}^{p-2}\big)  \nn \\
&\le& C_\eps\,N^\eps  \sup_{|a_n|\le1} \Big\|\summe a_n\xi_n(\omega_2) e(n\cdot) \Big\|_{2(p-2)}^{p-2}\;\|g-h\|_q,\nn \\
&=:& C_\eps\,N^\eps  J_{2(p-2),N}^{\half}(\omega_2)\;\|g-h\|_q,\label{eq:J2(p-2)}
\end{eqnarray}
where the $N^{\eps}$-term follows from Bernstein's inequality and we have set
\[ 
\sup_{|a_n|\le1} \Big\|\summe a_n\xi_n(\omega_2) e(n\cdot) \Big\|_{2(p-2)}^{2(p-2)} =: J_{2(p-2),N}(\omega_2).
\]
As before, $g,h\in \sqrt{m}\calP_A$, $A=A(\omega_2)=\{n\in[1,N]\:|\: \xi_n(\omega_2)=1\}$, 
see~\eqref{eq:PAdef}. One concludes that, for $\eps>0$ small and $q<\infty$ large depending on $\eps$, 
\bea
\log N_2(\calE(\omega_2),t) &\le& \log N_q\Big(\calP_{A(\omega_2)}, N^{-\eps} m^{-\half}\,J_{2(p-2),N}^{-\half}\,t\Big) \nn \\
&\le& C_q\,m 
\left\{ \begin{array}{ll} 1 + \log \frac{1}{t} & \text{\ \ if\ \ }0<t< N^\eps \sqrt{m\,J_{2(p-2),N}(\omega_2)} \\
                                    (m^{-\half}\,J_{2(p-2),N}^{-\half}(\omega_2)\, N^{-\eps} t)^{-\nu} & 
 \text{\ \ if\ \ } t> N^\eps  \sqrt{m\,J_{2(p-2),N}(\omega_2)} 
                 \end{array}
         \right.  \nn
\end{eqnarray}
where $\nu>2$, see Corollary~\ref{cor:small_t} and Corollary~\ref{cor:entropy2}.  
Inserting this estimate into the last term of~\eqref{eq:step4_p} yields by the random
majorant property on~$2(p-2)\ge 2$, 
\begin{equation}
\label{eq:en_ab}
\Erw'_{\omega_2} \int_0^\infty \sqrt{\log N_2(\calE(\omega_2),t)}\,dt 
\le C_\eps\, N^\eps \,m \sqrt{\Erw I_{2(p-2),N}}
\end{equation}
and therefore finally, by Lemma~\ref{lem:norm1},
\bea
&&\Erw_\omega \sup_{|a_n|\le 1} \int_0^1 \Big| \summe a_n\xi_n(\omega) e(n\theta) \Big|^p \, d\theta \nn\\
&& \le C_\eps\, N^\eps \Big[ m^{\frac{p}{2}}+m^{\frac12} \Erw I_{p-1,N} + (1+\sqrt{m\tau}) \sqrt{\Erw I_{2(p-1),N}} + m \sqrt{\Erw I_{2(p-2),N}} \;\Big]\nn \\
&& \le C_\eps\, N^\eps \Big[ (\tau N)^{\frac{p}{2}}+(\tau N)^{\frac12} \Big(\tau^{p-1}N^{p-2}+(\tau N)^{\frac{p-1}{2}}\Big) \nn\\ 
&& \quad + (1+\tau \sqrt{N}) \Big(\tau^{2(p-1)}N^{2p-3}+(\tau N)^{p-1}\Big)^{\half} +
\tau N \Big(\tau^{2(p-2)}N^{2p-5}+ (\tau N)^{p-2}\Big)^{\half} \Big]\nn \\
&&  \le C_\eps\, N^\eps\Big[ \tau^p N^{p-1} +  \tau^{p-1}N^{p-\frac{3}{2}}+ (\tau N)^{\frac{p}{2}}\Big].
\label{eq:three_piece}
\end{eqnarray}
If $\tau\ge N^{-\half}$, then $\tau^pN^{p-1}\ge \tau^{p-1}N^{p-\frac32}$. Moreover, if
$\tau\le N^{\frac{3-p}{p-2}}$, then $\tau^{p-1}N^{p-\frac32}\le (\tau N)^{\frac{p}{2}}$. 
In particular, if $3\le p\le 4$, then $\tau^{p-1}N^{p-\frac{3}{2}} \les \Erw I_{p,N}$, 
and the result follows.
On the other hand, if $p\ge4$, then $\tau\ge N^{-\half}$ insures that $\tau^{p-1}N^{p-\frac{3}{2}} \les \tau^pN^{p-1}\les\Erw I_{p,N}$, as claimed. 

\noindent It remains to discuss $2\le p\le 3$. In that case, Lemma~\ref{lem:decoup} implies that
\bea
2^{-p} \int_0^1 \left|\summe a_n \xi_n e(n\theta) \right|^p \,d\theta &=& 
\int \int_0^1 \sum_{n\in R^1_t} a_n\xi_n e(n\theta) \sum_{k\in R^2_t} \bar a_k\xi_k e(-k\theta) 
\Big|\sum_{\ell\in R^2_t} a_\ell\xi_\ell e(\ell\theta) \Big|^{p-2}\, d\theta\,dt \nn \\
&& + O\left( m^{\frac{p}{2}} \int_0^1 \Big(1+\Big|\summe \frac{a_n}{\sqrt{m}}\xi_n e(n\theta)\Bigr|^{2} \Big)\, d\theta \right). \label{eq:step1_p''}
\end{eqnarray}
The integral 
in~\eqref{eq:step1_p''} is~$O(1)$. Hence~\eqref{eq:step2_p} changes to
\begin{equation}
\Erw_\omega \sup_{|a_n|\le 1} \int_0^1 \Big| \summe a_n\xi_n(\omega) e(n\theta) \Big|^p \, d\theta 
\les m^{\frac{p}{2}} 
 + \Erw'_{\omega_2} \Erw'_{\omega_1}  \sup_{x\in \calE(\omega_2)} \sup_{|A|=m} 
\sum_{n\in A} \xi_n(\omega_1) x_n, \label{eq:step2_p'}
\end{equation}
with the same $\calE(\omega_2)$, and \eqref{eq:step4_p} becomes
\bea
&&\Erw_\omega \sup_{|a_n|\le 1} \int_0^1 \Big| \summe a_n\xi_n(\omega) e(n\theta) \Big|^p \, d\theta \nn\\
&& \le C_\eps\,N^\eps\Big[m^{\frac{p}{2}} + (1+\sqrt{m\tau})  \sqrt{\Erw I_{2(p-1),N}} +
\Erw_{\omega_2}' \int_0^\infty \sqrt{\log N_2(\calE(\omega_2),t)}\,dt\Big]. 
\label{eq:step4_p'}
\end{eqnarray}
Finally, the entropy estimate simplifies as $2(p-2)\le 2$ in this case: If $g|g|^{p-2},h|h|^{p-2}\in\calE(\omega_2)$, then $g,h\in\calP_{A(\omega_2)}$ and thus
\bea
 \|g|g|^{p-2}-h|h|^{p-2}\|_2 &\les& \|g-h\|_\infty \big(\|g\|_{2(p-2)}^{p-2} + \|h\|_{2(p-2)}^{p-2}\big) \nn \\
&\le& C_\eps\,N^{{\eps}} \|g-h\|_q \big(\|g\|_{2}^{p-2} + \|h\|_{2}^{p-2}\big)  \nn \\
&\le& C_\eps\,N^\eps\, m^{\frac{p-2}{2}} \|g-h\|_q, \nn
\end{eqnarray}
so that now
\[
\Erw'_{\omega_2} \int_0^\infty \sqrt{\log N_2(\calE(\omega_2),t)}\,dt 
\le C_\eps\, N^\eps \,m^{\frac{p}{2}}.
\]
We leave it to the reader to check that this again leads to~\eqref{eq:three_piece}.
As already mentioned above, the term $\tau^{p-1} N^{p-\frac32}$ can be absorbed into
$(\tau N)^{\frac{p}{2}}$, since $p\le 3$. 
\end{proof}

\noindent This lemma quickly leads to a proof of Theorem~\ref{thm:Lp} in case $\delta\le\half$
for $p>4$, and for all $0<\delta<1$ if $2<p<4$.  

\begin{cor}
\label{cor:delta<half}
Suppose $0<\delta\le \half$ and assume otherwise that the hypotheses of Theorem~\ref{thm:Lp}
are satisfied. Then~\eqref{eq:rmp} holds for all~$p\ge 4$. 
If $2<p<4$, then~\eqref{eq:rmp} holds for all~$0<\delta<1$. 
In particular, Theorem~\ref{thm:Lp} is valid in these cases.
\end{cor}
\begin{proof} As a first step, note that Lemma~\ref{lem:main_step} 
immediately implies that all odd integers
satisfy~\eqref{eq:rmp}. Next, one checks that~\eqref{eq:rmp} holds at~$p=\frac52$
since~$2(p-1)=3$ in that case. Now Lemma~\ref{lem:main_step}  implies that~\eqref{eq:rmp}
holds at all other values $p=\frac{2\ell+1}{2}$, for all integers $\ell\ge 3$. 
Generally speaking, one checks by means of induction that~\eqref{eq:rmp}
holds at all 
\[ p \in \Big\{ 2+\frac{\ell}{2^{j}}\:\Big|\: \ell \in \Z^+ \Big\}=:\calP_j. \]
Indeed, we just verified that this holds for $j=0,1$. Now assume that it holds
up to some integer~$j$ and we will prove it for~$j+1$. Thus take $p=2+\frac{\ell}{2^{j+1}}\in\calP_{j+1}$
such that~$2<p<3$. Then $2(p-1)=2+\frac{\ell}{2^{j}}$ for which~\eqref{eq:rmp} holds
by assumption. Hence Lemma~\ref{lem:main_step} applies. Now suppose $p\in\calP_{j+1}$
is such that $3<p<4$. Then~\eqref{eq:rmp} holds at~$p-1$ by what we just did, and at~$2(p-1),2(p-2)$
by assumption. Hence Lemma~\ref{lem:main_step} applies again. One now continues with $4<p<5$ etc.,
and we are done. Given any~$\eps>0$ and $p>2$ one can find $p_1<p<p_2$ with $p_1,p_2\in\calP_{j}$
where $p_2-p_1<\eps$. Hence~\eqref{eq:rmp} holds for all~$p$ by interpolation, as desired.
It remains to deal with $\delta>\half$ if $2<p<4$. Fix such a~$p$. Then by Bourgain's theorem on random~$\Lambda(p)$
sets, $\delta>\half$ implies that the random set~$S(\omega)$ is a~$\Lambda(p)$ set.
More precisely, 
\[ \Erw  \sup_{|a_n|\le 1}  \Big\| \summe a_n\xi_n(\omega) e(n\cdot) \Big\|_p^p \les (\tau N)^{\frac{p}{2}},\]
see Theorem~\ref{thm:jean} below.
Clearly,
\[  \Big\| \summe \xi_n(\omega) e(n\cdot) \Big\|_p \ge \Big\| \summe \xi_n(\omega) e(n\cdot) \Big\|_2 =
\#(S(\omega))^{\half},\]
and we are done.
\end{proof}

\subsection{Some improvements and  $\delta>\half$}

It is clear that the proof of Lemma~\ref{lem:main_step} in its present form does not allow us
to deal with the case $\delta>\half$. The difficulty arises from
the use of Plancherel in~\eqref{eq:L2kmaj} and~\eqref{eq:J2(p-2)}. Indeed,
once the $L^2$ bound is used, the estimates in the proof of Lemma~\ref{lem:main_step}
are optimal and they produce the unwanted~$\tau^{p-1}N^{p-\frac32}$ term in~\eqref{eq:three_piece}.
In order to improve this step, one can invoke Bourgain's theorem
on random~$\Lambda(p)$ sets. Recall the main theorem from~\cite{B1}:

\begin{theorem}
\label{thm:jean}
Fix some $p>2$. 
Let $\{\xi_j\}_{j=1}^N$ be selector variables as in Theorem~\ref{thm:Lp} 
with $\delta= 1-\frac{2}{p}$. Define
\begin{equation}
\label{eq:Kp}
 K_p(\omega)=\sup_{|a|_{\ell^2_N}\le1} \Big\|\summe a_n\xi_n(\omega) e(n\cdot)\Bigr\|_p
\end{equation}
Then $\Erw K_p^p \le C_p <\infty$. 
\end{theorem}

\noindent Although the main theorem in~\cite{B1} is formulated for generic sets rather
than in terms of expected values, this statement appears implicitly in~\cite{B1}, see 
page~241 (especially the last line on that page), as well as Section~5 of that paper. We will need
the following dual version of~\eqref{eq:Kp}. With $S(\omega)=\{n\in[1,N]\:|\: \xi_n(\omega)=1\}$,
\begin{equation}
\label{eq:dual_Kp}
\Big( \sum_{n\in S(\omega)} |\hat{f}(n)|^2 \Big)^{\half} \le K_p(\omega)\,\|f\|_{L^{p'}(\tor)}
\end{equation}
where $\frac{1}{p}+\frac{1}{p'}=1$. 

\begin{prop}
\label{prop:L5} 
If $4<p\le 7$, 
then the random majorant property~\eqref{eq:rmp} 
and therefore~\eqref{eq:prob} hold for all $0<\delta<1$.
\end{prop}
\begin{proof} It suffices to consider $\delta>\half$. 
This will be done in several steps.
For the sake of clarity, we first present the case $p=5$, and then indicate how to pass to the 
range~$4<p\le\frac{11}{2}$. We will then refine the argument even further to obtain the specified range.
The idea is to factor through a $\Lambda(3)$ set, i.e., 
in order to generate a random subset $S(\omega)\subset[1,N]$
of cardinality roughly~$N^{1-\delta}$ one first chooses a random subset $S_0(\omega)\subset[1,N]$
of cardinality about~$N^{\frac23}$, and then generates $S(\omega)\subset S_0(\omega)$. 
Hence, we let $\xi_j=\xi_j^{(0)}\,\xi^{(1)}_j$ where $\Erw \xi_j^{(0)}=N^{-\frac13}$, and 
$\Erw \xi_j^{(1)}=\tau'$ so that $\tau=N^{-\frac13}\tau'$. Moreover, we of course choose all these
random variables to be independent. The set $S_0(\omega_0):= \{n\in[1,N] \:|\: \xi_j^{(0)}(\omega_0)=1\}$
satisfies~\eqref{eq:Kp} and also its dual version~\eqref{eq:Kp} at~$p=3$. The argument is
similar to those in Theorems~\ref{thm:L3} and Lemma~\ref{lem:main_step}, so we will only
indicate those places that are different. Starting the argument as before, one arrives at
\begin{equation}
\Erw_\omega \sup_{|a_n|\le 1} \int_0^1 \Big| \summe a_n\xi_n(\omega) e(n\theta) \Big|^5 \, d\theta 
\les m^{\frac{5}{2}} +  m^{\half}  \Erw I_{4,N} 
 + \Erw'_{\omega_2} \Erw_{\omega_0}\Erw_{\omega_1}'  \sup_{x\in \calE(\omega_0,\omega_2)} \sup_{|A|=m} 
\sum_{n\in A} \xi_n^{(1)}(\omega_1)\, x_n, \label{eq:step2_5}
\end{equation}
in place of~\eqref{eq:step2_p}, where now
\[ 
\calE(\omega_0,\omega_2) = \Bigl\{ \Bigl(\Bigl|\Bigl\la e(n\cdot), 
\sum_{k=1}^N \bar b_k\xi_k(\omega_2) e(-k\cdot)\Big|\sum_{\ell=1}^N b_\ell\xi_\ell(\omega_2) e(\ell\cdot) \Big|^{3} \Bigr \ra \Bigr|\,\chi_{S_0(\omega_0)}(n) \Bigr)_{n=1}^N \:\Big|\: \sup_{1\le n\le N}|b_n|\le 1\Bigr\} 
\subset \R^N_+.
\]
Using \eqref{eq:dual_Kp} with $p'=\frac32$ instead of Plancherel and the majorant property at $p=6$ 
leads to
\bea
 \sup_{x\in\calE(\omega_0,\omega_2)} |x|_{\ell^2_N} &\le& K_3(\omega_0)\,\sup_{|a_k|\le1}  
\Big\|\sum_{k} a_k\,\xi_k (\omega_2) e(k\theta)\Bigr\|_6^4 \nn \\
&\le&  K_3(\omega_0)\, \Big\|\sum_{k} \xi_k (\omega_2) e(k\theta)\Bigr\|_6^4 = K_3(\omega_0)\, 
\Big(I_{6,N}(\omega_2)\Big)^{\frac23}, \nn
\end{eqnarray}
and thus
\begin{equation}
\label{eq:B_5}
 \Erw \sup_{x\in\calE(\omega_0,\omega_2)} |x|_{\ell^2_N} \le C\, \Big(\Erw I_{6,N}\Big)^{\frac23}.
\end{equation}
In the next step we use Lemma~\ref{lem:suplem} to bound the last term in~\eqref{eq:step2_5}
 for fixed $\omega_0,\omega_2$. Since $\sqrt{m\tau'}=N^{\frac16}\sqrt{m\tau}$, one obtains
from that lemma that 
\begin{equation}
\label{eq:sup5}
 \Erw_{\omega_1} \sup_{x\in \calE(\omega_0,\omega_2)} \sup_{|A|=m} \sum_{n\in A} \xi^{(1)}_n(\omega_1) x_n 
\les (1+N^{\frac16}\sqrt{m\tau})\sup_{x\in\calE(\omega_0,\omega_2)} |x|_{\ell^2_N} 
+\int_0^\infty \sqrt{\log N_2(\calE(\omega_0,\omega_2),t)}\,dt. 
\end{equation}
To control the entropy term, observe that 
the distance between any two elements in $\calE(\omega_0,\omega_2)$ is bounded by
\bea
 K_3(\omega_0)\,\|g|g|^{3}-h|h|^{3}\|_{\frac32} &\les& K_3(\omega_0)\,\|g-h\|_\infty \big(\|g\|_{\frac92}^{3} + 
\|h\|_{\frac92}^{3}\big) \nn \\
&\le& K_3(\omega_0)\,C_\eps\,N^{\eps} \|g-h\|_q \big(\|g\|_{5}^{3} + \|h\|_{5}^{3}\big)  \nn \\
&\le& K_3(\omega_0)\,C_\eps\,N^\eps  \sup_{|a_n|\le1} \Big\|\summe a_n\xi_n(\omega_2) e(n\cdot) \Big\|_{5}^{3}\;\|g-h\|_q. \nn 
\end{eqnarray}
As before, $q$ is large depending on $\eps>0$, $g,h\in \sqrt{m}\calP_A$, $A=A(\omega_2)=\{n\in[1,N]\:|\: \xi_n(\omega_2)=1\}$. 
Using Corollaries~\ref{cor:small_t} and~\ref{cor:entropy2}, one now arrives at
\bea
&& \Erw_{\omega_0}\,\Erw'_{\omega_2} \int_0^\infty \sqrt{\log N_2(\calE(\omega_0,\omega_2),t)}\,dt \nn \\ 
&& \le \Big(\Erw_{\omega_0} K_3(\omega_0)^3\Big)^{\frac13}\, C_\eps\, N^\eps \,m 
\left(\Erw_{\omega_2}\sup_{|a_n|\le1} \Big\|\summe a_n\xi_n(\omega_2) e(n\cdot) \Big\|_{5}^{5}\right)^{\frac35}. \label{eq:en_5}
\end{eqnarray}
Combining \eqref{eq:step2_5}, \eqref{eq:sup5}, \eqref{eq:B_5}, \eqref{eq:en_5} one obtains
\bea
\Erw_\omega \sup_{|a_n|\le 1}  \Big\| \summe a_n\xi_n(\omega) e(n\theta) \Big\|_5^5 
&\les&  (\tau N)^{\frac{5}{2}} +  (\tau N)^{\half} (\tau^4 N^3 + \tau^2 N^2) \nn \\
&&  + (1+N^{\frac16}\sqrt{N\tau^2})(\tau^6 N^5 + \tau^3 N^3)^{\frac23} \nn \\
&&  + 
C_\eps\, N^\eps \,N\tau 
\left(\Erw_{\omega_2}\sup_{|a_n|\le1} \Big\|\summe a_n\xi_n(\omega_2) e(n\cdot) \Big\|_{5}^{5}\right)^{\frac35} \nn \\
&\les&  C_\eps\, N^\eps (\tau N)^{\frac{5}{2}} + \tau^5 N^4 + \tau^{\frac92}N^{\frac72} + \tau^4N^{\frac{10}{3}} + N^{\frac83}\tau^3 \nn \\
&& + \half \Erw_\omega \sup_{|a_n|\le 1}  \Big\| \summe a_n\xi_n(\omega) e(n\theta) \Big\|_5^5. \nn
\end{eqnarray}
We leave it to the reader to check that the expressions with fractional exponents are dominated
by~$\tau^5 N^4$ provided~$\tau>N^{-\frac23}$. Hence, for those~$\tau$,  we have proved
\[ \Erw_\omega \sup_{|a_n|\le 1}  \Big\| \summe a_n\xi_n(\omega) e(n\theta) \Big\|_5^5 \le 
C_\eps\, N^\eps \Erw I_{5,N},
\]
and thus the RMP at $p=5$ and therefore also~\eqref{eq:prob} with~$p=5$ holds for all $0<\delta<\frac23$. 
On the other hand, if $\tau\le N^{-\frac23}<N^{-\frac35}$,
then $S(\omega)$ is a~$\Lambda(5)$ set with large probability by Bourgain's theorem. More precisely,
\eqref{eq:Kp} holds with~$p=5$ and thus~\eqref{eq:prob} follows with~$p=5$ for the range~$\delta\ge \frac23$ as well.

\noindent We now indicate how to obtain the range $4<p<\frac{11}{2}$. Instead of factoring through
a $\Lambda(3)$-set, one factors through a $\Lambda(q)$-set where $q'(p-1)=6$. Since we need
to cover the range $\half<\delta<1$, one needs to allow all $N^{-1}<\tau<N^{-\half}$. On the
other hand, the factorization means that $\tau=N^{-1+\frac{2}{q}}\tau'$ with some $\tau'<1$. 
This implies that necessarily $N^{-1+\frac{2}{q}}\ge N^{-\half}$ or $2\le q\le 4$. 
Hence $q'\ge\frac43$, and thus $\frac43(p-1)\le q'(p-1)= 6\Impl p\le\frac{11}{2}$. Inspection of the previous argument reveals that
we also need $q'(p-2)<p$, which by our choice of~$q'$ is the same as $p^2-7p+12 = (p-4)(p-3)>0$. 
But this holds for all $p>4$. We are now ready to run the same argument as before.
Observe that the first step already requires the (random) majorant property at~$p-1$. Therefore,
we start with the range $4<p<5$ so that this property is ensured by Corollary~\ref{cor:delta<half}. 
Analogously to the case~$p=5$ one arrives at the bound
\bea
M_p^p:=\Erw_\omega \sup_{|a_n|\le 1}  \Big\| \summe a_n\xi_n(\omega) e(n\theta) \Big\|_p^p 
&\les&  (\tau N)^{\frac{p}{2}} +  C_\eps\,N^\eps(\tau N)^{\half} \Big(\tau^{p-1} N^{p-2} + (\tau N)^{\frac{p-1}{2}}\Big) \nn \\
&&  + \Big(1+N^{\frac{1}{q'}}\tau\Big)(\tau^6 N^5 + \tau^3 N^3)^{\frac{p-1}{6}} 
  + C_\eps\, N^\eps \,N\tau \;M_p^{p-2} \nn \\
&\les&  C_\eps\, N^\eps\Big[ (\tau N)^{\frac{p}{2}} + \tau^p N^{p-1}\Big] + \tau^{p-1}N^{\frac56(p-1)} \nn \\
&& + 
\tau^{\frac{p+1}{2}}N^{\frac{p-1}{2}+\frac{1}{q'}} + \half M_p^p. \nn 
\end{eqnarray}
One now checks that the two unwanted terms in the final expression are dominated by $\tau^pN^{p-1}$
provided $\tau>N^{-\frac{p-1}{6}}$ and $\tau>N^{-\frac23}$, respectively. By the usual reduction
to the (random) $\Lambda(p)$-property, it suffices to consider the range $\tau>N^{-1+\frac{2}{p}}$. 
But since for $p\ge4$
\[ 1-\frac{2}{p}\le\frac{p-1}{6} \Longleftrightarrow p\ge 4,\qquad 1-\frac{2}{p} \le \frac23 
\Longleftrightarrow p\le6,\]
we are done with the case $4<p\le 5$. Finally, if $5<p<\frac{11}{2}$, then we just showed that
the RMP holds at $p-1$, and so we can repeat the exact same argument.

\noindent Next, we increase $p$ even further.  
For example, take $p=7$ and factor through a~$\Lambda(3)$ set. More precisely, set $q'=\frac43$. Then
$q'(p-1)=8$ and $q'(p-2)<p$. Since the majorant property holds at~$p-1=6$, one can repeat the 
argument for $p=5$  mutatis mutandis (use $L^8$ instead of~$L^6$). 
We leave it to the reader to check that this leads to
\[ 
M_7^7 \les C_\eps N^\eps\Big[\tau^7 N^6 + (\tau N)^{\frac72}\Big] + \tau^6 N^{\frac{21}{4}} 
+ \tau^4 N^{\frac{15}{4}}. 
\]
Moreover, the unwanted terms are dominated by $\tau^7 N^6$ provided $\tau>N^{-\frac34}$. But
since $1-\frac{2}{p}=\frac57 < \frac34$, the remaining range of~$\tau$ is covered by the 
random~$\Lambda(7)$  property as before.
Using~$L^8$ instead of~$L^6$ of course gives a larger range of~$p$'s. 
Indeed, let now $\frac43 \le q'\le 2$ be general such that $q'(p-1)=8$. 
This is possible for all $5\le p\le 7$. 
On the other hand, we also require $q'(p-2)<p$. This reduces
to $0<p^2-9p+16$ which means that $p>\half(9+\sqrt{17})$. Finally, to run the argument
we also need to know the RMP at~$p-1$. This was clear for $p=7$,
but it is not if $p$~is below~$7$. However, we will show in the next step that the RMP at~$p=7$
allows us to increase the range of $p$ from~$\frac{11}{2}$ (which would be
insufficient for our purposes) to~$\frac{25}{4}>6$. This in turn settles the issue
of $p-1$ if~$p<7$. Thus one does indeed obtain the RMP for all
$ \half(9+\sqrt{17}) \le p \le 7$.

\noindent Next, we argue that the RMP at $p=7$ makes it possible to increase~$p$ 
from~$\frac{11}{2}$ to~$\frac{25}{4}$.  To be precise, set $q'(p-1)=7$. The restriction
$\frac43\le q'\le 2$ yields~$\frac{9}{2}\le p\le \frac{25}{4}$. On the other hand,
$q'(p-2)<p$ is the same as $0<p^2-8p+14$, which holds for all $p>4+\sqrt{2}$ and thus,
in particular, for $p\ge\frac{11}{2}$. Finally, if $p\le\frac{25}{4}$, then
 the random majorant property~\eqref{eq:rmp}
holds at $p-1\le\frac{21}{4}<\frac{11}{2}$ by the first part of the proof. 
We can now run the same argument as before to conclude that
\bea
M_p^p
&\les&  (\tau N)^{\frac{p}{2}} +  C_\eps\,N^\eps(\tau N)^{\half} \Big(\tau^{p-1} N^{p-2} + (\tau N)^{\frac{p-1}{2}}\Big) \nn \\
&&  + C_\eps\,N^\eps\, \Big(1+N^{\frac{1}{q'}}\tau\Big)\Big(\tau^7 N^6 + (\tau N)^{\frac72}\Big)^{\frac{p-1}{7}} 
  + C_\eps\, N^\eps \,N\tau \;M_p^{p-2} \nn \\
&\les&  C_\eps\, N^\eps\Big[ (\tau N)^{\frac{p}{2}} + \tau^p N^{p-1} + \tau^{p-1}N^{\frac67(p-1)}  + 
\tau^{\frac{p+1}{2}}N^{\frac{p-1}{2}+\frac{1}{q'}}\Big] + \half M_p^p. \nn 
\end{eqnarray}
One now checks that the two unwanted terms in the final expression are dominated by $\tau^pN^{p-1}$
provided $\tau>N^{-\frac{p-1}{7}}$ and $\tau>N^{-\frac57}$, respectively. By the usual reduction
to the (random) $\Lambda(p)$-property, it suffices to consider the range $\tau>N^{-1+\frac{2}{p}}$. 
But since for $p\ge4$
\[ 1-\frac{2}{p}\le\frac{p-1}{7} \Longleftrightarrow p\ge 4+\sqrt{2},\qquad 1-\frac{2}{p} \le \frac57
\Longleftrightarrow p\le7,\]
we are done with this case as well.

\noindent It remains to close
 the gap $\frac{25}{4} < p < \half(9+\sqrt{17})$. The idea is to forfeit the 
requirement $q'(p-2)<p$ and instead replace
$q'(p-2)$ by the smallest number~$r$ to the right of~$q'(p-2)$ for which the random majorant
property is known. One then  uses H\"older's inequality which brings in $\Erw I_{r,N}$ instead of~$M_p$.
In our case the best choice of $r$ is~$r=\half(9+\sqrt{17})$. 
To be precise, we  set $q'(p-1)=8$ which by $q'\le\frac43$ can be done for $p\le 7$. But we are
only interested in $p\le r$, which is equivalent to $q'(p-2)=8\frac{p-2}{p-1}\le r$. Since we
already know that the RMP holds at~$p-1$, we can proceed as before, but
using H\"older's inequality to pass from~$q'(p-2)$ to~$r$. One checks that this leads to
\bea
M_p^p
&\les&  (\tau N)^{\frac{p}{2}} + C_\eps\,N^\eps\,(\tau N)^{\half} \Big(\tau^{p-1} N^{p-2} + (\tau N)^{\frac{p-1}{2}}\Big) 
 +  \Big(1+N^{\frac{1}{q'}}\tau\Big)\Big(\tau^8 N^7 + (\tau N)^{4}\Big)^{\frac{p-1}{8}} \nn \\
&&   +  N\tau \Big(\tau^r N^{r-1}+(\tau N)^{\frac{r}{2}}\Big)^{\frac{p-2}{r}}  \nn \\
&\les&  (\tau N)^{\frac{p}{2}} + \tau^p N^{p-1} + \tau^{p-1}N^{\frac78(p-1)} +N^{\frac{1}{q'}}\tau(\tau N)^{\frac{p-1}{2}}   + \tau^{p-1} N^{\frac{p-1}{r'}}N^{\frac{1}{r}}. \nn
\end{eqnarray}
This yields the desired bound under the conditions $\tau>N^{-\frac34}$, $\tau>N^{-\frac{p-1}{8}}$,
$\tau>N^{-\frac{p-2}{r}}$. Since $\frac34\ge\frac{p-1}{8}\ge\frac{p-2}{r}$, this reduces to
$\delta\le \frac{p-2}{r}$. To recapitulate, for the range $\frac{25}{4}\le p\le r$ we have
raised the admissible values of~$\delta$ from~$\half$ to~$\frac{p-2}{r}$, which is at 
least~$\frac{17}{4r}$. The point is now 
that this allows us to factor through~$\Lambda(q)$-sets for values of~$q$ larger than~$4$. Indeed,
define $q_0$ by $1-\frac{2}{q_0}=\frac{17}{4r}$. Going back to the argument involving the random majorant 
property at~$p=7$, we see that we can apply it for all $p$ for which $q'(p-1)=7$ with~$2\le q\le q_0$.
Recall that $q'(p-2)<p$ holds if $p>4+\sqrt{2}$, and thus for all~$p$ in the range under consideration. 
In order to close the gap $[\frac{25}{4},r]$ we therefore only need to check that
$q_0'(r-1)\le 7$. But since 
 $\frac{1}{q_0'}=\frac12+\frac{17}{8r}$, this is the same as $8r^2-36r-119\le0$. One explicitly
checks that with $r=\half(9+\sqrt{17})$ one has $8r^2-36r-119<-10$, and we are done.
\end{proof}

\subsection{Choosing subsets by means of correlated selectors}

To conclude this section, we want to address the issue of obtaining a version of Theorem~\ref{thm:Lp}
for subsets which are obtained by means of selectors $\xi_j$ that 
are allowed to have some degree of dependence. More precisely,
we will work with the selectors from the following definition.

\begin{defi}
\label{def:double}
Let $0<\tau<1$ be fixed. 
Define~$\xi_j(\omega)=\chi_{[0,\tau]}(2^j\omega)$
for $j\ge 1$. Here $\omega\in\tor=\R/\Z$ with probability measure $\Prob(d\omega)=d\omega$
equal to normalized Lebesgue measure. 
\end{defi}

Since the doubling map $\omega\mapsto 2\omega\;\mod 1$
is measure preserving, it follows that $\Erw \xi_j=\tau$ and $\Prob[\xi=1]=\tau$, $\Prob[\xi_j=0]=1-\tau$,
as in the random case. However, these selector variables are no longer 
independent. Nevertheless, they are close enough to being independent to make the following theorem 
accessible to the methods of the previous section.

\begin{theorem}
\label{thm:double}
Let $0<\delta<1$ be fixed. For every positive integer $N$ we let $\xi_j=\chi_{[0,\tau]}(2^j\omega)$ 
be as in Definition~\ref{def:double}  with $\tau=N^{-\delta}$. Define a subset
\begin{equation}
\label{eq:S_doub}
 S(\omega)=\{j\in[1,N]\:|\: \xi_j(\omega)=1\}
\end{equation}
for every $\omega\in\tor$. Then for every $\eps>0$ and $7\ge p\ge2$  one has 
\begin{equation}
\label{eq:prob_doub} 
\Prob\Bigl[\sup_{|a_n|\le1} \Big\|\sum_{n\in S(\omega)} a_n e(n\theta) \Big\|_{L^p(\tor)} \ge N^\eps\,
             \Big\|\sum_{n\in S(\omega)}  e(n\theta) \Big\|_{L^p(\tor)} \Bigr] \to 0
\end{equation}
as $N\to\infty$. Moreover, under the additional 
restriction~$\delta\le\half$, \eqref{eq:prob_doub} holds for all $p\ge 7$.
\end{theorem}

To prove this theorem we may of course assume that $\tau=2^{-k}$ for some positive integer~$k$. 
Then $\xi_j$ is measurable with respect to the dyadic intervals of length~$2^{-k-j}$ on
the unit interval~$\tor$, denoted by~$\calD_{j+k}$. Moreover, it is easy to see that
$\xi_{j}$ and~$\xi_{j+ak}$ are independent variables. 

\begin{lemma}
\label{lem:doub_indep}
Fix $j\ge0$ and $k\ge1$. Let $\tau=2^{-k}$ and $\xi_i$ be as in Definition~\ref{def:double}.
Then the sequence $\Big\{\xi_{j+ak}\Big\}_{a=1}^\infty$ is a realization
of a $0,1$-valued Bernoulli sequence with $\Erw\xi_i=\tau$. 
\end{lemma}
\begin{proof}
Fix $a>1$ and note that the variable $\xi_{j+ak}(\omega)$ is $2^{-(j+ak)}$-periodic.
On the other hand, each of the variables $\xi_{j+bk}$ with $b<a$ is constant on intervals
from~$\calD_{j+ak}$ (which is the same as saying that these variables are all~$\calD_{j+ak}$
 measurable). It follows that 
\[ \Prob\big[\xi_{j+ak} = 1\:|\: \xi_{j+bk}=\eps_b,\;0\le b\le a-1\big] = \tau
= \Prob[\xi_{j+ak} = 1],\]
for any choice of $\eps_b=0,1$, $0\le b\le a-1$. This implies independence.
\end{proof}

From now on, let $\tau=N^{-\delta}$ for some fixed $0<\delta<1$. 
In view of Lemma~\ref{lem:doub_indep} we can decompose the sequence $\{\xi_j\}_{j=1}^N$ 
into about~$\log N$ many subsequences, where the indices run along arithmetic 
progressions $\calP_i$ of step-size equal to~$\sim \log N$, and $1\le i\les \log N$. 
Each of the subsequences consists of i.i.d.~variables, but variables from different 
subsequences are not independent. This easily shows that Lemma~\ref{lem:norm1}
remains valid here, possibly with a logarithmic loss in the upper bound for~$\Erw I_{p,N}$.
Indeed, recall that the proof of that Lemma is based upon splitting a random trigonometric
polynomial into its expectation and a mean-zero part. Since the $L^p$-norm of the
Dirichlet kernel on an arithmetic progression of length~$K$ is about~$K^{\frac{1}{p'}}$, 
and here $\#\calP_i \sim \frac{N}{\log N}$, one sees immediately that the upper bound 
from~\eqref{eq:obere} is the same up to logarithmic factors. As far as the lower bound of
Lemma~\ref{lem:norm1} is concerned, note that the proof relies on obtaining {\em upper}
bounds on certain error terms, cf.~\eqref{eq:untere}-\eqref{eq:LDT2}. However, these
upper bounds are again immediate corollaries of the random case by virtue of the splitting
into the progressions~$\calP_i$. 

The consequence of this is that basically all the main estimates from the previous section
remain valid here, up to possibly an extra factor of~$\log N$. Clearly, such factors
are irrelevant in this context. More precisely, with $\xi_j$ as in Definition~\ref{def:double}
and~$S(\omega)$ as in~\eqref{eq:S_doub}, it is a corollary of the proof of Theorem~\ref{thm:Lp} that
\begin{equation}
\label{eq:up_doub}
\Erw \sup_{|a_n|\le1} \Big\|\sum_{n\in S(\omega)} a_n e(n\theta) \Big\|^p_{L^p(\tor)} 
\le C_\eps\, N^{\eps} \big(\tau^p N^{p-1}+(\tau N)^{\frac{p}{2}}\big).
\end{equation}
The proof of Theorem~\ref{thm:double} is therefore  completed
as before by appealing to (the adapted version) of Lemma~\ref{lem:norm1}.

\begin{remark} 
\label{rem:powers}
Other examples of much more strongly correlated selectors  
are $\xi_j(\omega)=\chi_{[0,\tau]}(j^s\omega)$ where $s$ is a fixed
positive integer and $\omega\in\tor$. It appears to be rather difficult to prove
a version of Theorem~\ref{thm:Lp} for these types of selectors.
\end{remark}

\section{Perturbing arithmetic progressions}
\label{sec:RMP2}

Let $\calP\subset [1,N]$ be an arithmetic progression of length~$L$, i.e., 
\[\calP=\{b+a\ell \:|\: 0<b<a,\;0\le \ell < L:=\lfloor N/a\rfloor \}\subset [1,N]. \]
Fix some arbitrary $\eps_0>0$. Suppose $N^{\eps_0} < s < a$ and 
let $\{\xi_j\}_{j\in\calP}$ be i.i.d.~variables, integer valued and uniformly 
distributed in~$[-s,s]$. We define a random subset 
\begin{equation}
\label{eq:arith_per}
\calS(\omega) := \{ j+\xi_j(\omega)\:|\: j\in\calP\}.
\end{equation}
For future reference, we set $I_j:= [j-s,j+s]$ for each $j\in\calP$. By construction, 
$S(\omega)\subset \bigcup_{j\in\calP} I_j$, and the intervals $I_j$ are congruent and
pairwise disjoint.

\subsection{Suprema of random processes}

The following lemma is related to Lemma~\ref{lem:suplem}.

\begin{lemma}
\label{lem:suplem_2} 
Let $\calE\subset \R_+^N$, $B=\sup_{x\in\calE} |x|$, 
and $\calS(\omega)$ be as in~\eqref{eq:arith_per}. Then
\[
\Erw_\omega \sup_{x\in\calE}  \sum_{j\in \calS(\omega)}  x_j \les
B\big(1+\sqrt{{L/s}}\big) + \int_0^B \sqrt{\log N_2(\calE,t)}\, dt
\]
where $N_2$ refers to the $L^2$ entropy.
\end{lemma}
\begin{proof} 
As in the proof of Lemma~\ref{lem:suplem}, we introduce $2^{-k}$-nets $\calE_k$ and~$\calF_k\subset\R^N$
so that $\diam(\calF_k)\le 1$,
\begin{equation}
\label{eq:F_kcard_2}
\log\,\#\calF_k \le C\,\log\,\#\calE_{k+1},
\end{equation}
and 
\begin{equation}
\label{eq:reduc_2}
\Erw \sup_{x\in\calE} \sum_{n\in \calS(\omega)}  x_n
\le \sum_{k\ge k_0} 2^{-k+1}  \Erw \sup_{y\in\calF_k} \;\sum_{n\in\calS(\omega)}  |y_n|.
\end{equation}
Now fix some $k\ge k_0$ and write $\calF$ instead of $\calF_k$. 
With $0<\rho_2$ to be determined, one has for any $|y|\le1$
\[
\sum_{i\in \calS(\omega)}  y_i \le \sum_{y_i\ge \rho_2} y_i + \sum_{y_i<\rho_2} \chi_{\calS(\omega)}(i)\, y_i 
\le \rho_2^{-1}  + \sum_{y_i<\rho_2} \chi_{\calS(\omega)}(i)\, y_i. 
\]
Let $q:=1+\lfloor\log\calF\rfloor$. Then, as in~\eqref{eq:sup_F},
\begin{equation}
\Erw\sup_{y\in\calF}\sum_{i\in \calS(\omega)}  y_i 
\les \rho_2^{-1} + \sup_{|y|\le 1}\Bigl\|\sum_{y_i<\rho_2} \chi_{\calS(\omega)}(i)\, y_i\Bigr\|_{L^q(\omega)} \label{eq:sup_F_2}.
\end{equation}
To control the last term in~\eqref{eq:sup_F_2}, we need the following analogue of~\eqref{eq:q_sum}.
By the multinomial theorem (for any positive integer~$q$),
\bea
&& \Erw \Big[\sum_{n\in\calS(\omega)} \chi_A(n) \Big]^q = \Erw \Bigl [ \sum_{j\in\calP} \chi_A(j+\xi_j(\omega)) \Bigr]^q = 
\sum_{q_1+\ldots+q_L=q} \binom{q}{q_1,\ldots,q_L} \Erw \prod_{j\in\calP} \chi_A(j+\xi_j(\omega))^{q_j}
\nn \\
&=&  \sum_{\nu=1}^q \;\sum_{\substack{1\le i_1<i_2<\ldots < i_\nu\le L\\ i_1,\ldots,i_\nu\in\calP}} \;
\sum_{\substack{q_{i_1}+\ldots+q_{i_\nu}=q\\ q_{i_1}\ge 1,\ldots,q_{i_\nu}\ge 1}} \binom{q}{q_{i_1},\ldots,q_{i_\nu}} \Erw \prod_{t=1}^\nu\chi_A(i_t+\xi_{i_t}(\omega)) \nn \\
&=&  \sum_{\nu=1}^q \;\;\sum_{\substack{1\le i_1<i_2<\ldots < i_\nu\le L\\ i_1,\ldots,i_\nu\in\calP}} \;
\sum_{\substack{q_{i_1}+\ldots+q_{i_\nu}=q\\ q_{i_1}\ge 1,\ldots,q_{i_\nu}\ge 1}} \binom{q}{q_{i_1},\ldots,q_{i_\nu}} \;\frac{|A\cap I_{i_1}|}{|I_{i_1}|} \frac{|A\cap I_{i_2}|}{|I_{i_2}|}\cdot\ldots\cdot \frac{|A\cap I_{i_\nu}|}{|I_{i_\nu}|} \nn \\  
&\le& \sum_{\nu=1}^q \nu^q \frac{1}{\nu!} \Bigl(\sum_{j\in\calP} \frac{|A\cap I_j|}{|I_j|} \Big)^\nu
\le \sum_{\nu=1}^q \frac{\nu^q}{q!} \frac{q!}{\nu!} \Bigl( \frac{|A\cap \bigcup_{j\in\calP}I_j|}{2s+1} \Big)^\nu \nn\\
&\le& \sum_{\nu=1}^q \binom{q}{\nu} q^{q-\nu} \Bigl( \frac{e|A\cap \bigcup_{j\in\calP}I_j|}{2s+1} \Big)^\nu
\le \Big(q+\frac{e|A\cap \bigcup_{j\in\calP}I_j|}{2s+1} \Big)^q. \nn
\end{eqnarray}
Continuing with the final term in \eqref{eq:sup_F_2} one concludes that
\bea
 \sup_{|y|\le 1}\Bigl\|\sum_{y_i<\rho_2} \chi_{\calS(\omega)}(i)\, y_i\Bigr\|_{L^q(\omega)} &\les&
\sum_{\rho_2^{-2}<2^j} 2^{-\frac{j}{2}} \sup_{|A|=2^j}\Bigl\| \sum_{n\in \calS(\omega)}\chi_A(n)  \Bigr\|_{L^q(\omega)} \nn \\
&\les& \sum_{\rho_2^{-2}<2^j} 2^{-\frac{j}{2}} \Big(q+\frac{\min(Ls,2^j)}{s}\Big) \les q\rho_2 + \sqrt{L/s}. \nn
\end{eqnarray}
Let~$\rho_2=q^{-\half}=(1+\log\#\calF)^{-\half}$.  
Inserting this bound into~\eqref{eq:sup_F_2} therefore yields
\[
\Erw\sup_{y\in\calF}\sum_{i\in \calS(\omega)}  y_i \les  \sqrt{q} + \sqrt{L/s} 
\les \sqrt{L/s} + 1 + \sqrt{\log\#\calF}. 
\]
The lemma now follows in view of \eqref{eq:F_kcard_2} and~\eqref{eq:reduc_2}.
\end{proof}

\subsection{The $L^p$ norm of the Dirichlet kernel over $\calS(\omega)$}

The following lemma determines an upper bound on the typical size of the 
Dirichlet kernel over~$\calS(\omega)$ in the $L^p$-norm,
with $2\le p\le 4$. The lower bound, as well as the case~$p>4$ 
will be dealt with below.

\begin{lemma}
\label{lem:norm2_upper}
With $\calS(\omega)$ as in~\eqref{eq:arith_per}, there exists a constant $C_p$ so that
\[ \Erw \Bigl\|\sum_{n\in\calS(\omega)} e(n\cdot) \Bigr\|_p^p 
\le C_p \Big(L^{\frac{p}{2}}+\frac{L^{p-1}}{s}\Big)
\]
for all $2\le p\le 4$. 
\end{lemma}
\begin{proof} 
For every $\ell\in\Z$ define
\[ 
A_\ell(\omega) := \#\{n,m\in\calS(\omega)\:|\: n-m=\ell \} = \sum_{j,k\in\calP} \chi_{[j-k+\xi_j-\xi_k=\ell]}.
\]
Clearly, $\calP-\calP \subset \bigcup_{i} J_i$ where $i\in a\Z$ and $J_i:=[i-2s,i+2s]$. These intervals
are mutually disjoint since $s\ll a$. This means that 
\[ 
\ell\in J_i \Impl A_\ell(\omega) = \sum_{j\in\calP} \chi_{[j-i\in\calP]} \chi_{[\xi_j-\xi_{j-i}=\ell-i]}. 
\]
Let us denote the unique $i$ for which $\ell\in J_i$ by $i(\ell)$. For simplicity, we shall mostly
write~$i$. 
If $i=0$, then $A_\ell(\omega)=L\delta_0(\ell)$ (recall that $\#\calP=L$). 
Otherwise, if $i\ne0$, then one finds that
\bea
\label{eq:ErwAell}
\Erw A_\ell &=& \sum_{j\in\calP} \frac{2}{2s+1} \Big(1-\frac{|\ell-i|}{s}\Big)_{+} \chi_\calP(j-i)
=   (L-|i|/a)_{+}\, \frac{2}{2s+1} \Big(1-\frac{|\ell-i|}{s}\Big)_{+} \\
&=& \frac{2L}{2s+1}\,\what{K_s}(\ell-i(\ell)) \what{K_L}(|i|/a) \label{eq:fejer}
\end{eqnarray}
where $\what{K_n}(k) = (1-|k|/n)_+$ denotes the Fejer kernel.
Moreover, if $i\ne0$, then 
\bea
\Erw A_\ell^2 &=& \Erw \sum_{\substack{j,k\in\calP\\j-i\in\calP,\,k-i\in\calP}} \chi_{[\xi_j-\xi_{j-i}=\ell-i]} \chi_{[\xi_k-\xi_{k-i}=\ell-i]} \nn \\
&=& \sum_{\substack{j,k\in\calP\\j-i\in\calP,\,k-i\in\calP}} \chi_{[j\ne k,j\ne k\pm i]} \Erw \chi_{[\xi_j-\xi_{j-i}=\ell-i]} \Erw\chi_{[\xi_k-\xi_{k-i}=\ell-i]}  \nn \\
&& + \sum_{\substack{j,k\in\calP\\j-i\in\calP,\,k-i\in\calP}} \Big(\chi_{[j=k,j\ne k\pm i]} +
\chi_{[j\ne k, k+ i,j=k-i]} +\chi_{[j\ne k, k-i,j= k+ i]}\Big) \Erw \chi_{[\xi_j-\xi_{j-i}=\ell-i]} \chi_{[\xi_k-\xi_{k-i}=\ell-i]}. \nn 
\end{eqnarray}
Hence 
\bea
\Erw A_\ell^2 &=& \sum_{\substack{j,k\in\calP\\j-i\in\calP,\,k-i\in\calP}}  \Erw \chi_{[\xi_j-\xi_{j-i}=\ell-i]} \Erw\chi_{[\xi_k-\xi_{k-i}=\ell-i]}  \nn \\
&& + \sum_{\substack{j,k\in\calP\\j-i\in\calP,\,k-i\in\calP}} \Big(\chi_{[j=k,j\ne k\pm i]} +
\chi_{[j\ne k, k+ i,j=k-i]} +\chi_{[j\ne k, k-i,j= k+ i]}\Big) \Erw\Big(\chi_{[\xi_j-\xi_{j-i}=\ell-i]} \chi_{[\xi_k-\xi_{k-i}=\ell-i]}\Big) \label{eq:coinc1} \\
&& - \sum_{\substack{j,k\in\calP\\j-i\in\calP,\,k-i\in\calP}} \Big(\chi_{[j=k,j\ne k\pm i]} +
\chi_{[j\ne k, k+ i,j=k-i]} +\chi_{[j\ne k, k-i,j= k+ i]}\Big) \Erw \chi_{[\xi_j-\xi_{j-i}=\ell-i]} \Erw \chi_{[\xi_k-\xi_{k-i}=\ell-i]} \label{eq:coinc2} \\
&& \quad = (\Erw  A_\ell)^2 + O\Big(\frac{L}{s}\Big(1-\frac{|\ell-i|}{s}\Big)_+\Big). \label{eq:erwAell2}
\end{eqnarray}
The $O$-term in \eqref{eq:erwAell2} arises because the error terms in~\eqref{eq:coinc1} and~\eqref{eq:coinc2} basically reduce to the computation of a single expectation as in~\eqref{eq:ErwAell}. Now consider
\[ 
V_{p,N}:= \int_0^1 \Big|\sum_{\ell\in\Z} (A_\ell(\omega)-\Erw A_\ell) e(\ell\theta) \Big|^{\frac{p}{2}}
\,d\theta. 
\]
Since $p\le 4$ by assumption, $\Erw V_{p,N} \le (\Erw V_{4,N})^{\frac{p}{4}}$. Moreover,
by~\eqref{eq:erwAell2}, 
\bea
\Erw V_{4,N} &=& \Erw \sum_{\ell\in\Z} |A_\ell(\omega)-\Erw A_\ell|^2 = \sum_{\ell\in\Z} \Big[\Erw(A_\ell^2) - (\Erw A_\ell)^2 \Big] \nn \\
&=& \Erw (A_0^2) - (\Erw A_0)^2 + \sum_{\ell\ne0}\left[ (\Erw A_\ell)^2 + O\Big(\frac{L}{s}\Big(1-\frac{|\ell-i|}{2s+1}\Big)_+\Big)\right] - \sum_{\ell\ne0} (\Erw A_\ell)^2 \nn \\
&\les& L^2 \nn
\end{eqnarray}
and therefore 
\begin{equation}
\label{eq:erwVpN}
\Erw V_{p,N} \les L^{\frac{p}{2}}.
\end{equation}
 In view of~\eqref{eq:fejer},
\bea
\sum_{\ell\in\Z} \Erw A_\ell\; e(\ell\theta) &=& \sum_{\ell\in\Z} \frac{2L}{2s+1} \,\what{K_s}(\ell-i(\ell)) \what{K_L}(|i(\ell)|/a) e((\ell-i(\ell))\theta) e(i(\ell)\theta) \nn \\
&=& \frac{2L}{2s+1} \sum_{k\in\Z} \what{K_s}(k) e(k\theta) \sum_{j\in \Z} \what{K_L}(j) e(ja\theta) 
= \frac{2L}{2s+1} K_s(\theta) K_L(a\theta). \nn
\end{eqnarray}
It follows that
\bea
\int_0^1 \Big|\sum_{\ell\in\Z} \Erw A_\ell\;e(\ell\theta) \Big|^{\frac{p}{2}} \, d\theta
&\les& \Big(\frac{L}{s}\Big)^{\frac{p}{2}} \int_0^1 \Big(\frac{1}{s}\min(s^2,\theta^{-2})\Bigr)^{\frac{p}{2}} \; |K_L(a\theta)|^{\frac{p}{2}}\, d\theta \nn \\
&\les& \Big(\frac{L}{s}\Big)^{\frac{p}{2}} \Bigl\{s^{\frac{p}{2}}\,a^{-1}\,L^{\frac{p}{2}-1} + \sum_{j=1}^a\Big(\frac{1}{s}\min(s^2,(j/a)^{-2})\Bigr)^{\frac{p}{2}}\,a^{-1}\,L^{\frac{p}{2}-1} \Bigr\} \nn \\
&\les& \frac{L^{p-1}}{s}. \label{eq:sumerw}
\end{eqnarray}
Combining~\eqref{eq:erwVpN} with \eqref{eq:sumerw} one obtains for $2\le p\le 4$
\bea
&& \Erw \int_0^1 \Bigl|\sum_{n\in\calS(\omega)} e(n\theta) \Bigr|^p \,d\theta = \Erw 
\int_0^1 \Bigl|\sum_{\ell\in\Z} A_\ell(\omega) e(\ell\theta) \Bigr|^{\frac{p}{2}} \,d\theta \nn \\
&&  \les
\int_0^1 \Bigl|\sum_{\ell\in\Z} \Erw A_\ell\; e(\ell\theta) \Bigr|^{\frac{p}{2}} \,d\theta 
+ \Erw \int_0^1 \Bigl|\sum_{\ell\in\Z} [A_\ell(\omega)-\Erw A_\ell] e(\ell\theta) \Bigr|^{\frac{p}{2}} \,d\theta \nn \\
&& \les \frac{L^{p-1}}{s} + L^{\frac{p}{2}}, \label{eq:Lp-1}
\end{eqnarray}
as claimed. 
\end{proof}

The following lemma is a special case of a well-known large deviation
estimate for martingales with bounded increments. The norm $\|\cdot\|_\infty$ 
refers to the supremum norm with respect to the probability space.

\begin{lemma}
\label{lem:LinftyLDT} 
Suppose $\{X_j\}_{j=1}^M$ are complex-valued independent variables with $\Erw X_j=0$. 
Then for all $\lambda>0$ 
\[ \Prob\Big[ \big|\sum_{j=1}^M X_j \big| > \lambda \Bigl(\sum_{j=1}^M \|X_j\|^2_\infty\Big)^{\half} \Big]
< C\,e^{-c\lambda^2}
\]
with some absolute constants $c,C$. 
\end{lemma}

Lemma~\ref{lem:LinftyLDT} implies the following simple generalization of the Salem-Zygmund bound.

\begin{cor}
\label{cor:SZ2}
Let $s,L$ be positive integers. 
Suppose $T_L$ is a trigonometric polynomial with random coefficients that can be written in the form
\[ T_L(\theta) = \sum_{j=-L}^L a_j(\theta) e(j\theta)\]
where $a_j(\theta)$ are trigonometric polynomials of degree at most~$s$,
and such that for fixed~$\theta$ they are independent random variables with $\Erw a_j(\theta)=0$. Moreover,
we assume that $\sup_{\theta\in\tor}|a_j(\theta)|\le1$ for each~$j$. Then for every $A>1$ 
\[ \Prob[ \|T_L\|_\infty > C\sqrt{\log (s+L)}\, \sqrt{L}] \le (s+L)^{-A},\]
with some constant $C=C(A)$.
\end{cor}
\begin{proof} 
Fix $\theta\in\tor$ and apply Lemma~\ref{lem:LinftyLDT} with $X_j=a_j(\theta)e(j\theta)$. 
By assumption, these are complex valued independent mean-zero variables with $\|X_j\|_\infty\le1$.
Therefore,
\begin{equation}
\label{eq:TNfixtheta}
 \sup_{\theta\in\tor}\; \Prob\Big[ \big|\sum_{j=-L}^L a_j(\theta)e(j\theta) \big| > \lambda \sqrt{L} \Big] 
< C\,e^{-c\lambda^2}.
\end{equation}
If $|\theta-\theta'|<(s+L)^{-2}$, then by Bernstein's inequality 
\[|T_L(\theta)-T_L(\theta')|\le (s+L)\|T_L\|_\infty |\theta-\theta'|\les (s+L)L(s+L)^{-2} \les 1.\]
Now pick a $(s+L)^{-2}$-net on the circle. The corollary follows by setting $\lambda=C\log(s+L)$
with~$C$ large, and summing~\eqref{eq:TNfixtheta} over the elements of the net.
\end{proof}

We can now state the general version of Lemma~\ref{lem:norm2_upper}. 
It is possible to remove the $\log$-term from the upper bound,
but the bound given below suffices for our purposes.

\begin{lemma}
\label{lem:norm2_ap}
For all $p\ge2$ there exists $C_p$ so that
\begin{equation}
\label{eq:upper_loss} 
\Erw \Bigl\|\sum_{n\in\calS(\omega)} e(n\cdot) \Bigr\|_p^p 
\le C_p \Big(\frac{L^{p-1}}{s}+(L\,\log N)^{\frac{p}{2}}\Big).
\end{equation}
Moreover, there is $c_p>0$ small so that
\[ 
\Prob\Big[ \Big\|\sum_{n\in\calS(\omega)} e(n\cdot) \Bigr\|_p^p < c_p 
\Big(L^{\frac{p}{2}}+\frac{L^{p-1}}{s}\Big) \Big]
\to 0\]
as $N\to\infty$. 
\end{lemma}
\begin{proof}
We work with the following splitting:
\begin{equation}
\label{eq:alter_split}
\sum_{n\in\calS(\omega)} e(n\theta) = \sum_{n\in\Z} \Erw \chi_{\calS(\omega)}(n) e(n\theta)
+ \sum_{n\in\Z} \Big[ \chi_{\calS(\omega)}(n)  - \Erw \chi_{\calS(\omega)}(n) \Big] e(n\theta).
\end{equation}
Clearly,
\begin{equation} 
\label{eq:dir_again}
\sum_{n\in\Z} \Erw \chi_{\calS(\omega)}(n)e(n\theta) = 
\frac{1}{2s+1} D_s(\theta) \sum_{j\in\calP} e(j\theta),
\end{equation}
and thus
\bea
\Big\| \sum_{n\in\Z} \Erw \chi_{\calS(\omega)}(n)e(n\theta) \Big \|_p^p 
&\les& s^{-p} \int_0^1 \Big| \min(s,\theta^{-1})\, \sum_{j=1}^L e(ja\theta) \Big|^p\,d\theta \nn \\
&\les& s^{-p} \Big[ \sum_{k=1}^L \min(s,a/k)^{p} + s^p\,  \Big]\frac{L^{p-1}}{a} 
\les \frac{L^{p-1}}{s}. \label{eq:dir_upper}
\end{eqnarray}
Conversely, 
\bea
\Big\| \sum_{n\in\Z} \Erw \chi_{\calS(\omega)}(n)e(n\theta) \Big \|_p^p 
&\gtr& s^{-p} \int_0^{1/s} \Big| D_s(\theta) \sum_{j=1}^L e(ja\theta) \Big|^p \, d\theta \nn \\
&\gtr& \frac{a}{s} \frac{L^{p-1}}{a} = \frac{L^{p-1}}{s}. 
\label{eq:dir_lower}
\end{eqnarray}
Both \eqref{eq:dir_upper} and~\eqref{eq:dir_lower} hold for all $p>1$. 
The second sum in~\eqref{eq:alter_split} can be written as 
\[
\sum_{n\in\Z} \Big[ \chi_{\calS(\omega)}(n)  - \Erw \chi_{\calS(\omega)}(n) \Big] e(n\theta)
= \sum_{j\in\calP} a_j(\omega,\theta) e(j\theta),
\]
where $a_j(\omega,\theta) = \chi_{I_j}(\xi_j(\omega))e(\xi_j(\omega)) - \frac{1}{2s+1} D_s(\theta)$. Clearly, $\Erw a_j(\omega,\theta)=0$, 
$\sup_\theta |a_j(\omega,\theta)|\le 2$ and for fixed $\theta$ the random variables $a_j(\omega,\theta)$
are independent. Thus Corollary~\ref{cor:SZ2} yields that
\begin{equation}
\label{eq:mean0_infty}
\Big\|\sum_{n\in\Z} \Big[ \chi_{\calS(\omega)}(n)  - \Erw \chi_{\calS(\omega)}(n) \Big] e(n\theta)\Big\|_\infty \les \sqrt{L} \sqrt{\log N}
\end{equation}
up to probability at most $(s+L)^{-p}$. In particular,
\[
\Erw \Big\|\sum_{n\in\Z} \Big[ \chi_{\calS(\omega)}(n)  - \Erw \chi_{\calS(\omega)}(n) \Big] e(n\theta)\Big\|_p^p \les (L\,\log N)^{\frac{p}{2}} + L^p (s+L)^{-p} \les (L\,\log N)^{\frac{p}{2}}.
\]
In conjunction with~\eqref{eq:dir_upper} this yields~\eqref{eq:upper_loss}.
 For the lower bound, take $N^{-\eps_0/2}> h\gg\frac{1}{s}$. Then
\bea
&& \int_0^1 \Big|\sum_{n\in\calS(\omega)} e(n\theta) \Big|^p\,d\theta \gtr 
\int_0^{1/s} \Big|\sum_{n\in\Z} \Erw \chi_{\calS(\omega)}(n) e(n\theta)\Big|^p\,d\theta -
\int_0^{1/s} \Big|\sum_{n\in\Z} \Big[ \chi_{\calS(\omega)}(n)  - \Erw \chi_{\calS(\omega)}(n) \Big] e(n\theta)\Big|^p
\, d\theta \nn \\
&& \mbox{\hspace{1in}}+ \int_h^{1-h}  \Big|\sum_{n\in\Z} \Big[ \chi_{\calS(\omega)}(n)  - \Erw \chi_{\calS(\omega)}(n) \Big] e(n\theta)\Big|^p \, d\theta - \int_h^{1-h} \Big|\sum_{n\in\Z} \Erw \chi_{\calS(\omega)}(n) e(n\theta)\Big|^p\,d\theta \nn \\
&& \mbox{\hspace{1in}} =: I+II+III+IV.  \label{eq:I-IV}
\end{eqnarray}
By \eqref{eq:dir_lower}, $I\gtr \frac{L^{p-1}}{s}$. Secondly, 
\bea
IV &\les& \int_h^{1-h} \Big| \frac{1}{s}D_s(\theta) \sum_{j=0}^{L-1} e(ja\theta) \Big|^p\,d\theta \nn \\
&\les& s^{-p} \sum_{j>ah} (j/a)^{-p}\,\frac{L^{p-1}}{a} \les s^{-p} h^{-p+1} L^{p-1} \ll \frac{L^{p-1}}{s}
\label{eq:IV}
\end{eqnarray}
where the final estimate follows from $hs\gg1$. Thirdly, in view of $p\ge2$ and~\eqref{eq:mean0_infty},
\bea
III &\gtr& \int_0^1 \Big|\sum_{n\in\Z} \Big[ \chi_{\calS(\omega)}(n)  - \Erw \chi_{\calS(\omega)}(n) \Big] e(n\theta)\Big|^p \, d\theta \nn \\
&& \mbox{\hspace{1in}} - \Big(\int_0^h + \int_{1-h}^1\Big)\Big|\sum_{n\in\Z} \Big[ \chi_{\calS(\omega)}(n)  - \Erw \chi_{\calS(\omega)}(n) \Big] e(n\theta)\Big|^p \, d\theta \nn \\
&\gtr& \left(\int_0^1 \Big|\sum_{n\in\Z} \Big[ \chi_{\calS(\omega)}(n)  - \Erw \chi_{\calS(\omega)}(n) \Big] e(n\theta)\Big|^2 \, d\theta \right)^{\frac{p}{2}} - C\, h (L\,\log N)^{\frac{p}{2}} \nn \\
&\gtr& L^{\frac{p}{2}} - C\, h (L\,\log N)^{\frac{p}{2}}, \label{eq:III}
\end{eqnarray}
up to probability $(s+L)^{-p}=o(1)$ as $N\to\infty$. Similarly, \eqref{eq:mean0_infty} implies that
\[ II \les  s^{-1}\, (L\,\log N)^{\frac{p}{2}} \]
up to probability $(s+L)^{-p}$. 
Combining this bound with~\eqref{eq:III}, \eqref{eq:IV}, and~\eqref{eq:I-IV} implies that
\[ 
\int_0^1 \Big|\sum_{n\in\calS(\omega)} e(n\theta) \Big|^p\,d\theta \gtr \frac{L^{p-1}}{s} + L^{\frac{p}{2}}
- C(h+s^{-1}) (L\,\log N)^{\frac{p}{2}}
\]
asymptotically with probability one.
Since $h<N^{-\eps}$ and $s>N^\eps$,  the lemma follows.
\end{proof}

\subsection{The majorant property for randomly perturbed arithmetic progressions}

We are now ready to state our first result for perturbed arithmetic progressions as
defined in~\eqref{eq:arith_per}. In this section, if $\calS$ is the perturbation of 
an arithmetic progression of length~$L$, then we write 
\[ A_{p,L}(\omega) := \Big\|\sum_{n\in\calS(\omega)} e(n\cdot) \Bigr\|_p^p.\]
Also, we say that the random majorant property (RMP) holds at $p$ if
\begin{equation}
\label{eq:rmp_ap} 
\Erw_\omega \sup_{|a_n|\le1} \Big\|\sum_{n\in\calS(\omega)} a_n\, e(n\cdot) \Bigr\|_p^p \le C_\eps\,N^\eps 
\Erw_\omega  \Big\|\sum_{n\in\calS(\omega)}  e(n\theta) \Bigr\|_p^p.
\end{equation}
Of course, this depends on the length~$L$ of the underlying arithmetic progression.
Although~$L$ is arbitrary, it will be kept fixed in the course of any argument
that uses~\eqref{eq:rmp_ap}. 

\begin{theorem}
\label{thm:AP1}
Let $\calS$ be as in \eqref{eq:arith_per}.  
Then for every $\eps>0$ and $4\ge p\ge2$  one has 
\begin{equation}
\label{eq:prob_AP} 
\Prob\Bigl[\sup_{|a_n|\le1} \Big\|\sum_{n\in \calS(\omega)} a_n e(n\theta) \Big\|_{L^p(\tor)} \ge N^\eps\,
             \Big\|\sum_{n\in \calS(\omega)}  e(n\theta) \Big\|_{L^p(\tor)} \Bigr] \to 0
\end{equation}
as $N\to\infty$. Moreover, under the additional 
restriction~$L\ge s$, \eqref{eq:prob_AP} holds for all $p\ge 4$.
\end{theorem}
\begin{proof}
The proof is similar to the random case of the previous section, so
we shall be somewhat brief. We will show that the RMP holds at $p$ provided either 
$2\le p\le 3$, or if the RMP holds at $p-1$, $2(p-1)$, and~$2(p-2)$. 
It is important to notice that the RMP at~$p$ implies~\eqref{eq:prob_AP}.
Firstly, recall that we can write $\calS(\omega)=\{j+\xi_j\:|\:j\in\calP\}$. 
We apply the decoupling lemma, Lemma~\ref{lem:decoup}, to the progression~$\calP$.
I.e., in the notation of Lemma~\ref{lem:decoup}, $R^1_t:=\{j\in\calP \:|\:\zeta_j=1\}$,
and $R^2_t:=\{j\in\calP \:|\:\zeta_j=0\}$. Set
\[ \calS_t^1(\omega):=\{j+\xi_j(\omega)\:|\:j\in R^1_t\},\quad \calS_t^2(\omega):=\{j+\xi_j(\omega)\:|\:j\in R^2_t\}. \]  
Therefore, by Lemma~\ref{lem:decoup}, 
\bea
\frac18 \int_0^1 \Big|\sum_{n\in\calS(\omega)} a_n e(n\theta) \Big|^p \,d\theta &=& 
\int \int_0^1 \sum_{n\in \calS^1_t(\omega)} a_n\, e(n\theta) \sum_{k\in \calS^2_t(\omega)} \bar a_k\, e(-k\theta) \Big|\sum_{\ell\in \calS^2_t(\omega)} a_\ell\, e(\ell\theta) \Big|^{p-2}\, d\theta\,dt \nn \\
&& + O\left( L^{\frac{p}{2}} \int_0^1 \Big(1+\Big|\sum_{n\in\calS} \frac{a_n}{\sqrt{L}}\, e(n\theta)\Bigr|^{\max(p-1,2)} \Big)\, d\theta \right). \label{eq:step1_AP}
\end{eqnarray}
If either $p\le 3$, or if the RMP holds at $p-1$, then the $O$-term in~\eqref{eq:step1_AP} 
is at most
\begin{equation}
\label{eq:rmp_p-1}
 L^{\frac{p}{2}} + C_\eps\,N^{\eps}L^{\half} \Erw A_{p-1,L}  \les N^\eps L^{\frac{p}{2}},
\end{equation}
see Lemma~\ref{lem:norm2_ap}. We therefore obtain as in~\eqref{eq:step2}, 
\bea
&& \Erw_\omega \sup_{|a_n|\le 1} \int_0^1 \Big| \sum_{n\in\calS(\omega)} a_n\, e(n\theta) \Big|^p \, d\theta \nn \\
&\les& C_\eps\,N^{\eps}L^{\frac{p}{2}} 
+ \int \Erw_{\omega_1,\omega_2} \sup_{|a_n|\le 1} \left| \int_0^1  \sum_{n\in \calS^1_t(\omega_1)} a_n\, e(n\theta) \sum_{k\in \calS^2_t(\omega_2)} \bar a_k\, e(-k\theta) 
\;\Big|\sum_{\ell\in \calS^2_t(\omega_2)} a_\ell\, e(\ell\theta) \Big|^{p-2}\, d\theta \; \right|\,dt \nn \\
&\les& C_\eps\,N^{\eps}L^{\frac{p}{2}} 
+ \int \Erw_{\omega_1}\Erw_{\omega_2} \sup_{\substack{|a_n|\le 1\\|b_n|\le1}} \left| \int_0^1  
\sum_{n\in \calS^1_t(\omega_1)} a_n\, e(n\theta) \sum_{k\in \calS^2_t(\omega_2)} \bar b_k\, e(-k\theta) 
\;\Big|\sum_{\ell\in \calS^2_t(\omega_2)} b_\ell\, e(\ell\theta) \Big|^{p-2}\, d\theta \; \right|\,dt \nn \\
&\les& C_\eps\,N^{\eps}L^{\frac{p}{2}} + 
\int \Erw_{\omega_1}\Erw_{\omega_2} \sup_{\substack{|a_n|\le 1\\|b_n|\le1}}
\left|\int_0^1  \sum_{n\in\calS(\omega_1)} a_n\, e(n\theta) \sum_{k\in\calS(\omega_2)} \bar b_k\, e(-k\theta)\; \Big|\sum_{\ell\in\calS(\omega_2)} b_\ell\, e(\ell\theta) \Big|^{p-2}\, d\theta \; \right|\,dt \nn \\
&& \quad \les C_\eps\,N^{\eps}L^{\frac{p}{2}} + 
\Erw_{\omega_2} \Erw_{\omega_1}  \sup_{x\in \calE(\omega_2)} 
\sum_{n\in \calS(\omega_1)}  x_n.\label{eq:step2_AP}
\end{eqnarray}
Here
\[ 
\calE(\omega_2) := \Bigl\{ \Bigl(\Bigl|\Bigl\la e(n\cdot), 
\sum_{k\in\calS(\omega_2)} \bar b_k\, e(-k\cdot)\Big|\sum_{\ell\in\calS(\omega_2)} b_\ell\, e(\ell\cdot) \Big|^{p-2} \Bigr \ra \Bigr| \Bigr)_{n=1}^N \:\Big|\: \sup_{1\le n\le N}|b_n|\le 1\Bigr\} \subset \R^N_+.
\]
By Lemma~\ref{lem:suplem_2}, it follows from~\eqref{eq:step2_AP} that
\bea
\label{eq:meschdu}
&& \Erw_\omega \sup_{|a_n|\le 1} \int_0^1 \Big| \sum_{n\in\calS(\omega)} a_n\, e(n\theta) \Big|^p \, d\theta \\
&& \les C_\eps\,N^{\eps}L^{\frac{p}{2}} + (1+\sqrt{L/s})\Erw_{\omega_2}\sup_{x\in \calE(\omega_2)} |x|
+ \Erw_{\omega_2} \int_0^\infty \sqrt{\log N_2(\calE(\omega_2),t)}\,dt. \nn 
\end{eqnarray}
Now suppose the RMP holds at $2(p-1)$ (so this holds for sure if $p$ is an odd integer). Then by Plancherel, 
\[ \Erw_{\omega_2}\sup_{x\in \calE(\omega_2)} |x| \le C_\eps N^{\eps} \Erw_{\omega} \Big\| \sum_{n\in\calS(\omega)}  e(n\cdot) \Big\|^{p-1}_{2(p-1)} \le C_\eps N^{\eps} \sqrt{\Erw_{\omega} A_{2(p-1),L}(\omega)}.\]
As far as the entropy term in~\eqref{eq:meschdu} is concerned, the same analysis as in the random
case shows that if  $p\le3$, then
\[ \Erw_{\omega_2} \int_0^\infty \sqrt{\log N_2(\calE(\omega_2),t)}\,dt \le C_\eps\,N^\eps\, L^{\frac32},\]
or if $p>3$ and the RMP holds at $2(p-2)$, then
\[ \Erw_{\omega_2} \int_0^\infty \sqrt{\log N_2(\calE(\omega_2),t)}\,dt \le C_\eps\,N^\eps\, L\sqrt{\Erw A_{2(p-2),L}},\]
see \eqref{eq:J2(p-2)} and~\eqref{eq:en_ab} for the details. Inserting all of this into~\eqref{eq:meschdu}
yields, under the assumption that  $p>3$ and
the RMP holds at $p-1$, $2(p-1)$, and $2(p-2)$ (the case $p\le3$ is similar),  
\bea 
&&  \Erw_\omega \sup_{|a_n|\le 1} \int_0^1 \Big| \sum_{n\in\calS(\omega)} a_n\, e(n\theta) \Big|^p \, d\theta \nn \\
&\les& C_\eps\,N^{\eps}\left\{L^{\frac{p}{2}} + (1+\sqrt{L/s}) \sqrt{\Erw_{\omega} A_{2(p-1),L}(\omega)}
+  L\sqrt{\Erw A_{2(p-2),L}} \right\} \nn \\
&\les& C_\eps\,N^{\eps}\left\{L^{\frac{p}{2}} + (1+\sqrt{L/s}) \Big(\frac{L^{2p-3}}{s}+L^{p-1}\Big)^{\half}
+  L\Big(\frac{L^{2p-5}}{s}+L^{p-2}\Big)^{\half} \right\} \nn \\
&\les& C_\eps\,N^{\eps} \Big[ \frac{L^{p-1}}{s} + L^{\frac{p}{2}} + \frac{L^{p-\frac32}}{\sqrt{s}} \Big]. \label{eq:uff_ap}
\end{eqnarray}
Recall from Lemma~\ref{lem:norm2_ap} that the desired bound is $\frac{L^{p-1}}{s}+L^{\frac{p}{2}}$. 
If $p=3$, then~\eqref{eq:uff_ap} does indeed agree with this bound. Since the hypotheses involving the~RMP
hold in case~$p=3$, we are done with that case, regardless of the relative size of $L$ and~$s$. 
Let us assume now that $L\ge s$. Then~\eqref{eq:uff_ap} agrees with the desired bound for all~$p$.
This means that we can run the same type of inductive argument as in Corollary~\ref{cor:delta<half}.
We leave it to the reader to check that this proves~\eqref{eq:prob_AP} for all $p\ge2$ provided $L\ge s$.
Finally, if $L<s$, then $L<s\le a\le \frac{N}{L}$ and thus $L\le \sqrt{N}$. 
In particular, $\#\calS\le \sqrt{N}$ in that case. In analogy with the random subset case, this
suggests that~$\calS(\omega)$ are $\Lambda(p)$-sets for $2\le p\le 4$ with high probability.
Although perturbed arithmetic progressions are not covered by~\cite{B1}, it turns out that the
strategy from~\cite{B1} and~\cite{B2} is still relevant. More precisely, suppose first that $2\le p\le 3$.
Then~\eqref{eq:meschdu} holds, even without the~$N^\eps$-term. 
By Plancherel, but without appealing to any RMP,
\begin{equation}
 \Erw_{\omega_2}\sup_{x\in \calE(\omega_2)} |x| \le \Erw_{\omega_2} \sup_{|a_n|\le1} 
\Big\| \Big|\sum_{n\in\calS(\omega)} a_n e(n\theta) \Big|^{p-1} \Big\|_2 \les  K_p^{\frac{p}{2}} L^{\frac{p-2}{2}}. \label{eq:Kp_again}
\end{equation}
Here
\[ K_p^p := \Erw_\omega \sup_{|a_n|\le 1} \int_0^1 \Big| \sum_{n\in\calS(\omega)} a_n\, e(n\theta) \Big|^p \, d\theta. \]
To pass to \eqref{eq:Kp_again}, one writes $2(p-1)=p+(p-2)$ and then estimates 
the $(p-2)$-power in~$L^\infty$. Secondly, to bound the entropy term, set $q=\frac{2}{3-p}$. Then
by Plancherel the distance between any two elements in~$\calE(\omega_2)$ is at most
\bea
 \big\|g|g|^{p-2}-h|h|^{p-2}\big\|_2 &\les& \|g-h\|_q \big(\|g\|_{2}^{p-2} + \|h\|_{2}^{p-2}\big) \nn \\
&\les&  L^{\frac{p-2}{2}}\,\|g-h\|_q, \nn
\end{eqnarray}
where $g,h\in \sqrt{L}\calP_{\calS(\omega_2)}$, see~\eqref{eq:PAdef}. As before, the entropy estimate therefore reads
\[ 
\Erw_{\omega_2} \int_0^\infty \sqrt{\log N_2(\calE(\omega_2),t)}\,dt \les \sqrt{L}\, L^{\frac{p-2}{2}}\sqrt{L} = L^{\frac{p}{2}}.
\]
Inserting these bounds into~\eqref{eq:meschdu} yields
\[ K_p^p \les L^{\frac{p}{2}} + (1+\sqrt{L/s}) K_p^{\frac{p}{2}} L^{\frac{p-2}{2}} 
\le C L^{\frac{p}{2}} + \half K_p^p + C (1+L/s)L^{p-2}.\]
Since $L^{p-2}\le L^{\frac{p}{2}}$ in view of~$p\le 4$, one obtains the desired bound
\[ K_p^p \les L^{\frac{p}{2}} + \frac{L^{p-1}}{s}\]
if $2\le p \le3$ and regardless of the relative size of $L$ and $s$. 
If $3\le p\le 4$, then the previous argument needs to be modified in two places:
Firstly, there is the issue of the $O$-term in~\eqref{eq:step1_AP}. However, we just showed
that the RMP holds at~$p-1\le 3$, and therefore~\eqref{eq:rmp_p-1} applies here as well
 (even without the $N^\eps$-term). Secondly, the entropy bounds need to be modified.
In case $3\le p\le 4$, one has $2(p-2)\le p$. Hence 
\bea
 \|g|g|^{p-2}-h|h|^{p-2}\|_2 &\le& \|g-h\|_\infty \big(\|g\|_{2(p-2)}^{p-2} + \|h\|_{2(p-2)}^{p-2}\big) \nn \\
&\le& C_\eps\,N^{\eps} \|g-h\|_q \big(\|g\|_{p}^{p-2} + \|h\|_{p}^{p-2}\big)  \nn \\
&\le& C_\eps\,N^\eps  \sup_{|a_n|\le1} \Big\|\summe a_n\xi_n(\omega_2) e(n\cdot) \Big\|_{p}^{p-2}\;\|g-h\|_q,\nn
\end{eqnarray}
with $g,h$ as above. By the usual arguments, cf.~\eqref{eq:J2(p-2)}, it follows that
\[ 
\Erw_{\omega_2} \int_0^\infty \sqrt{\log N_2(\calE(\omega_2),t)}\,dt \les L\, K_p^{\frac{p-2}{2}} \le \half K_p^p + L^{\frac{p}{2}}.
\]
Inserting these bounds into \eqref{eq:meschdu} implies the desired bound.
\end{proof}

\begin{remark} It is possible that one can make improvements on Theorem~\ref{thm:AP1}
similar to those in Theorem~\ref{thm:Lp}, thus removing the condition~$L\ge s$ in some range
of $p\ge4$. This would require working with $\Lambda(p)$
 type arguments as we just did in the end of the previous proof.
But we do not pursue that issue here.
\end{remark}

\underline{Acknowledgement:} The second author was partially supported by an NSF grant, DMS-0070538,
and a Sloan fellowship.

\bibliographystyle{amsplain}

\begin{thebibliography}{99}


\bibitem[B1]{B1} Bourgain, J. {\it Bounded orthogonal systems and the $\Lambda(p)$-set problem.}  
Acta Math.\ 162  (1989),  no.\ 3-4, 227--245.

\bibitem[B2]{B2} Bourgain, J. {\it $\Lambda(p)$-sets in analysis: results, problems and related aspects. Handbook of the geometry of Banach spaces,}  Vol.~I, 195--232, North-Holland, Amsterdam, 2001. 

\bibitem[B3]{B3} Bourgain, J. {\it On $\Lambda(p)$-subsets of squares.}  Israel J. Math.  67  (1989),  no. 3, 291--311.

\bibitem[BLM]{BLM} Bourgain, J.; Lindenstrauss, J.; Milman, V. {\it Approximation of zonoids by zonotopes.}  Acta Math.\ 162  (1989),  no.\ 1-2, 73--141.

\bibitem[LT]{LT} Ledoux, M.; Talagrand, M. {\it Probability in Banach spaces. Isoperimetry and processes.}  Ergebnisse der Mathematik und ihrer Grenzgebiete (3) [Results in Mathematics and Related Areas (3)], 23. Springer-Verlag, Berlin, 1991. 

\bibitem[M]{M} Mockenhaupt, G. {\em Bounds in Lebesgue spaces of oscillatory integrals.} Habilitationsschrift, Siegen, 1996.
Available online, www.math.gatech.edu/\~{}gerdm/

\bibitem[PT-J]{PT} Pajor, A.; Tomczak-Jaegermann, N. {\it Subspaces of small codimension of finite-dimensional Banach spaces.}
  Proc.\ Amer.\ Math.\ Soc.\ 97  (1986),  no.~4, 637--642.

\bibitem[P]{P} Pisier, G. {\it The volume of convex bodies and Banach space geometry.}  Cambridge Tracts in Mathematics, 94. Cambridge University Press, Cambridge, 1989.

\bibitem[T]{T} Talagrand, M. {\em Sections of smooth convex bodies via majorizing measures.}   Acta Math.\ 175  (1995),  no.~2, 273--300.


\end{thebibliography}

\medskip\noindent
\textsc{Mockenhaupt: School of Mathematics, Georgia Tech, Atlanta, GA 30332, U.S.A.}
{\em email: }\textsf{\bf gerdm@math.gatech.edu} \\
\textsc{Schlag: Division of Astronomy, Mathematics, and Physics, 253-37 Caltech, Pasadena, CA 91125, U.S.A.}\\
{\em email: }\textsf{\bf schlag@its.caltech.edu}

\end{document}